\magnification=1200
\hoffset=.5pc
\voffset=.5pc
\input amssym.def
\input amssym.tex
 \font\newrm =cmr10 at 15pt
\def\bul{\raise .9pt\hbox{\newrm .\kern-.105em } }

\def\fr{\frak}

\baselineskip=13pt
 
 \def\h{\hbox{ }}
 
 \def\p{{\fr p}}
 \def\u{{\fr u}}

 \def\a{{\fr a}}
 \def\d{{\fr d}}
 
 \def\ss{{\fr s}}
 
 \def\b{{\fr b}}
 
 \def\hh{{\fr h}}

 \def\q{{\fr q}}

 \def\<{\le}
 \def\>{\ge}

 \def\s{{\h\subset\h}}
 
 \def\vs{\vskip }

 \def\mapright#1
  {\smash{\mathop
  {\longrightarrow}
  \limits^{#1}}}

 \def\kk#1{{\kern .4 em} #1}
 \def\vs{\vskip 1pc}

\hsize = 31pc
\vsize = 45pc
\overfullrule = 0pt

\overfullrule = 0pt

\font\ninerm=cmr9
\rm
\centerline{\newrm Gelfand--Zeitlin theory from the perspective}
\vskip 5pt
\centerline{\newrm  of classical mechanics  I}\vskip 1.1pc 
 \baselineskip=11pt
\vskip8pt
\centerline{\rm BERTRAM
KOSTANT\footnote*{\ninerm Research supported in part by NSF grant
DMS-0209473 and in part by the \hfil\break KG\&G Foundation.}}\vskip 1pc
\centerline{and}\vskip 1pc
\centerline{\rm NOLAN WALLACH\footnote{**}
{\ninerm
Research supported in part by NSF grant MTH 0200305.}}
\vskip 2pt
\footnote{}
{\ninerm \noindent 2000 Mathematics Subject Classification:
 Primary14L30, 14R20, 33C45, 53D17}\vskip 2pc\baselineskip=11pt

{\font\ninerm=cmr9 \baselineskip=14pt 
\noindent{\bf ABSTRACT.  }\ninerm \rm Let $M(n)$ be the algebra (both Lie and
associative) of $n\times n$ matrices over $\Bbb C$. Then $M(n)$ inherits a
Poisson structure from its dual using the bilinear form $(x,y) = -tr\,xy$. The $Gl(n)$
adjoint orbits are the symplectic leaves and the algebra, $P(n)$, of
polynomial functions on $M(n)$ is a Poisson algebra. In particular
if $f\in P(n)$ then there is a corresponding vector field $\xi_f$
on $M(n)$. If $m\leq n$ then $M(m)$ embeds as a Lie subalgebra
of $M(n)$ (upper left hand block) and $P(m)$ embeds as a Poisson subalgebra of
$P(n)$. Then, as an analogue of the Gelfand--Zeitlin algebra in the enveloping algebra of
$M(n)$, let $J(n)$ be the subalgebra of $P(n)$ generated by $P(m)^{Gl(m)}$ for
$m=1,\ldots,n$. One observes that
$$J(n) \cong P(1)^{Gl(1)}\otimes \cdots \otimes P(n)^{Gl(n)}$$ We prove that
$J(n)$ is a maximal Poisson commutative subalgebra of $P(n)$ and that for
any $p\in J(n)$ the holomorphic vector field $\xi_{f}$ is integrable and
generates a global one parameter group $\sigma_p(z)$ of holomorphic
transformations of $M(n)$. If $d(n)= n(n+1)/2$ then $J(n)$ is a polynomial
ring
$\Bbb C[p_1,\ldots,p_{d(n)}]$ and the vector fields $\xi_{p_i},\,i\in I_{d(n)}$,
span a commutative Lie algebra of dimension $d(n-1)$. Let $A$ be a corresponding
simply-connected Lie group so that $A\cong \Bbb C^{d(n-1)}$. Then $A$ operates on $M(n)$
by an action $\sigma$ so that if $a\in A$ then $$\sigma(a) = \sigma_{p_1}(z_1)\cdots
\sigma_{p_{d(n)}}(z_{d(n)})$$ where $a$ is the product of $exp\,\,z_i\,\xi_{p_i}$
for $i=1,\ldots,d(n)$.  We prove that the orbits of $A$ are independent of
the choice of the generators $p_i$. In addition we prove the following results
about this rather remarkable group action. (The
latter is a very extensive enlargement of an abelian group action
 introduced in [GS]). 

 (1) Let $x\in M(n)$. Then $A\cdot x$ is an orbit
of maximal dimension ($d(n-1)$) if and only if the differentials
$(dp_{i})_x,\,i=1,\ldots,d(n)$, are linearly independent. Such
$A$-orbits are explicitly determined.

(2) The orbits, $O_x$, of the adjoint action of $Gl(n)$ on $M(n)$ are $A$-stable and if
$O_x$ is an orbit of maximal dimension $(n(n-1))$, that is, if $x$ is regular, then the
$A$-orbits of dimension $d(n-1)$ in $O_x$ are the leaves of a polarization of a Zariski
open dense subset of the symplectic manifold $O_x$. 

 The results of the paper are related to the theory of orthogonal polynomials. Motivated
by the interlacing property of the zeros of neighboring orthogonal polynomials on $\Bbb
R$ we introduce a certain Zariski open subset $M_{\Omega}(n)$ of
$M(n)$ and prove 

(3) $M_{\Omega}(n)$ has the structure of $(\Bbb
C^{\times})^{d(n-1)}$ bundle over a $d(n)$-dimensional ``Lagrangian"
manifold. Moreover the fibers are maximal $A$-orbits.

In Part II of this two part paper we deal with an analogue of the
Gelfand--Kirillov conjecture. The fibration in (3) leads to
the construction of $n^2 +1$ functions (including a constant
function) in an algebraic extension of the function field of $M(n)$
which, under Poisson bracket, satisfies the commutation relations of
the direct sum of a 
$2\,d(n-1) +1$ dimensional Heisenberg Lie algebra and an
$n$-dimensional commutative Lie algebra.\vskip 1pc
\baselineskip 15pt

\rm

	\centerline{\bf 0. Introduction}\vskip 1pc 0.1. Let $M(n)$, for
any positive integer $n$, denote the Lie (and asssociative) algebra
of all complex $n\times n$ complex matrices. Let $P(n)$ be the
graded commutative algebra of all polynomial functions on $M(n)$.
The symmetric
algebra over $M(n)$, as one knows, is a Poisson
algebra. Using the bilinear form $(x,y) = -tr\,xy$ on $M(n)$ this
may be carried over to $P(n)$, defining on $P(n)$
the structure of a Poisson algebra and hence the structure of a Poisson manifold on $M(n)$. Consequently to
each
$p\in P(n)$  there is associated a holomorphic
vector field $\xi_p$ on $M(n)$ such that $$\xi_p\,q = [p,q]$$ where $q\in P(n)$
and
$[p,q]$ is Poisson bracket. 

For any positive integer $k$ put $d(k) = k(k+1)/2$ and let $I_k$ be the set
$\{1,\ldots,k\}$. If $m\in I_n$ we will regard $M(m)$ (upper
left hand $m\times m$ corner) as a Lie subalgebra of $M(n)$. As a ``classical mechanics" analogue to
the Gelfand--Zeitlin commutative subalgebra of the univeral enveloping algebra of $M(n)$,
let $J(n)$ be the subalgebra of $P(n)$ generated by $P(m)^{Gl(m)}$ for all $m\in I_n$. Then $$J(n) =
P(1)^{Gl(1)}\otimes
\cdots
\times P(n)^{Gl(n)}$$ In addition 
$J(n)$ is a Poisson commutative polynomial subalgebra of $P(n)$ with $d(n)$ generators. In fact we can write
$J(n) = \Bbb C[p_1,\ldots,p_{d(n)}]$ where, for $x\in M(n)$, $p_i(x), i\in I_{d(n)}$, 
``run over" the elementary symmetric functions of the roots of the characteristic polynomial of
$x_m,\,\,m\in I_n$. Here and throughout, $x_m\in M(m)$ is the upper left $m\times m$
minor of $x$. The algebraic morphism  $$\Phi_n: M(n) \to \Bbb C^{d(n)}\,\,\,\hbox{where}\,\,\Phi_n(x)
= (p_1(x),\ldots,p_{d(n)}(x))\eqno (0.1)$$ plays a major role in this paper. Let $\b_e$ be the
$d(n)$-dimensional affine space of all
$x\in M(n)$ of the form
$$x = 
\left(\matrix{a_{1\,1}&a_{1\,2}&\cdots &a_{1\,n-1}&a_{1\,n}\cr 1&a_{2\,2}&\cdots &a_{2\,n-1}&a_{2\,n}\cr
0&1&\cdots &a_{3\,n-1}&a_{3\,n}\cr
\vdots &\vdots &\ddots &\vdots &\vdots\cr 0&0&\cdots &1&a_{n\,n}\cr}\right)\eqno (0.2)$$ where $a_{i\,j}\in
\Bbb C$ are arbitrary. As a generalization of classical facts about companion matrices we prove
\vs {\bf Theorem 0.1.} {\it The restriction $$\b_e\to \Bbb C^{d(n)}\eqno (0.3)$$ of the map $\Phi_n$ is an
algebraic isomorphism.}
\vs See Theorem 2.3 and Remark 2.4. The real and imaginary parts of a complex number define a
lexicographical in $\Bbb C$. For any
$x\in M(n)$ and
$m\in I_n$ let
$E_x(m) =
\{\mu_{1\,m}(x),\ldots,\mu_{m\,m}(x)\}$ be the (increasing) ordered $m$-tuple of eigenvalues of 
$x_m$, with the multiplicity as roots of the characteristic polynomial. As a corollary of Theorem 0.1 one has
the following independence (with respect to $m$) of the eigenvalue sequences
$E_x(m)$.
\vs {\bf Theorem 0.2.} {\it For all
$m\in I_n$ let $E(m) =
\{\mu_{1\,m},\ldots,\mu_{m\,m}\}$ be an arbitrary $m$-tuple with values in $\Bbb C$. Then there exists a
unique
$x\in \b_e$ such that $E(m) = E_x(m)$, up to ordering, for all $m\in I_n$.}\vs See Theorem 2.5. For any $c\in
\Bbb C^{d(n)}$ let $M_c(n) = \Phi_n^{-1}(c)$ be the ``fiber" of $\Phi_n$ over $c$. If $x,y\in M(n)$ then 
$x$ and $y$ lie in the same fiber if and only if $E_x(m) = E_y(m)$ for all $m\in I_n$.
\vs {\bf Remark 0.3.} Theorem 0.1 implies that $\Phi_n$ is surjective and asserts that $\b_e$ is a
cross-section of 
$\Phi_n$. That is, to any
$c\in
\Bbb C^{d(n)}$ the intersection $M_c(n)\cap \b_e$ consists of exactly one matrix. \vs 0.2. The main results
of the present paper, Part I, of a two part paper, concern the properties of a complex analytic abelian group
$A$ of dimension
$d(n-1)$ which operates on
$M(n)$. One has \vs {\bf Theorem 0.4. } {\it The algebra $J(n)$ is a maximal Poisson commutative subalgebra of
$P(n)$. Furthermore the vector field $\xi_p$, for any $p\in J(n)$, is globally integrable on $M(n)$,
defining an analytic action of $\Bbb C$ on $M(n)$. Moreover the fiber $M_c(n)$ is stable under this action,
for any $c\in \Bbb C^{d(n)}$.} \vs See Theorem 3.4, Proposition 3.5 and Theorem 3.25. The generators $p_i$
of $J(n)$ are replaced by a more convenient set of $d(n)$ generators $p_{(i)},\,i\in I_{d(n)}$. See \S
3.1 and (3.20). The span of $\xi_{p_{(i)}},\,i\in I_{d(n)}$, is a commutative $d(n-1)$-dimensional Lie
algebra of analytic vector fields on $M(n)$. The Lie algebra $\a$ integrates to an action of a complex 
analytic group $A\cong \Bbb C^{d(n-1)}$ on $M(n)$. In a sense $A$ is very extensive enlargement of a group
introduced in
\S 4 of [GS]. However no diagonalizability, compactness or eigenvalue interlacing is required
for the existence of $A$. 

Now let $$x\mapsto V_x\eqno (0.4)$$ be the nonconstant dimensional tangent space ``distribution" on $M(n)$
defined by putting $V_x = \{(\xi_{p})_x\mid p\in J(n)\}$. Let $D(n)$ be the group of all diagonal matices
in $Gl(n)$ and let $Ad\,D(n)$ be the group of automorphisms of $M(n)$ defined by the adjoint action of
$D(n)$ on $M(n)$. \vs {\bf Theorem 0.5.} {\it The orbits of
$A$ are leaves of the distribution (0.4). Furthermore $M_c(n)$ is stable under the action of $A$, for
any
$c\in \Bbb C^{d(n)}$, so that
$M_c(n)$ is a union of $A$-orbits. Finally the abelian $d(n-1)$-dimensional group of analytic isomorphisms of
$M(n)$ (as an analytic manifold) defined by $A$, contains $Ad\,D(n)$ as an 
$n-1$-dimensional subgroup. }\vs See Theorem 3.4 and Proposition 3.5. For the last statement in the Theorem
0.5, see Theorem 3.27.
\vs For any
$x\in M(n)$ and
$m\in I_n$ let
$Z_{x,m}\s M(m)$ be the (obviously commutative) associative subalgebra generated by $x_m$ and the identity of
$M(m)$. Let
$G_{x,m}\s Gl(m)$ be the commutative algebraic subgroup of $Gl(n)$ corresponding to $Z_{x,m}$ when the latter
is regarded as a Lie algebra. The orbits of $A$ are described in \vs {\bf Theorem 0.6.} {\it Let
$x\in M(n)$. Consider the following morphism of nonsingular irreducible affine varieties
$$G_{x,1}\times
\cdots \times G_{x,n-1}
\to M(n)\eqno (0.5)$$ where for
$g(m)\in G_{x,m},\,\,m\in I_{n-1}$, $$ (g(1),\ldots, g(n-1))\mapsto Ad\,(g(1)\cdots g(n-1))(x)\eqno (0.6)$$
Then the image of (0.5) is exactly the
$A$-orbit $A\cdot x$. In particular
$A
\cdot x$ is an irreducible, constructible (in the sense of Chevalley) subset of $M(n)$. Furthermore the
 Zariski closure
of
$A\cdot x$ is the same as its closure in the Euclidean topology. In addition $A\cdot x$ contains a Zariski
open subset of its closure.} \vs See Theorem 3.6. One should note that even though the groups $G_{x,m},\,m\in
I_n$, are commutative (and $A$ is commutative) they do not commute with each other, so the order in (0.6) is
important. 

Let $x\in M(n)$. Then, by definition, $x$ is regular (in $M(n))$ if the $Gl(n)$ adjoint orbit of $x$
is of maximal dimension, $2\,d(n)$. One knows that $x$ is regular if and only if $dim\,Z_{x,n} = n$ or if
and only if the differentials $(dp_{d(n-1) + k})_x,\,k\in I_n$, are linearly independent where the indexing
of the $p_i$ is such that the ring of invariants, $P(m)^{Gl(m)}$, is generated by $p_{d(m-1) + k},\,k\in 
I_m$. We will now say that $x$ is strongly regular if $(dp_i)_x,\,\,i\in I_{d(n)}$, are linearly 
independent. 

One has that $dim\,A\cdot x\leq d(n-1)$. \vs {\bf Theorem 0.7.} {\it Let $x\in
M(n)$. Then the following conditions are equivalent. $$\eqalign{&(a)\,\,\,x\,\,\hbox{is strongly regular}\cr
&(b)\,\,\,A\cdot x\,\,\hbox{is an orbit of maximal dimension, $d(n-1)$}\cr&(c)\,\,\,dim\,Z_{x,m} =
m,\,\,\forall m\in I_n,\,\,\hbox{and  $Z_{x,m}\cap Z_{x,m+1}
= 0,\,\,\,\forall m\in I_{n-1}$}\cr}$$}\vs See Theorem 2.14 and (3.29). 

Let $M^{sreg}(n)\s M(n)$ be the Zariski open set of all strongly regular matrices. Note that $M^{sreg}(n)$
is not empty since in fact $\b_e\s M^{sreg}(n)$. Theorem 0.6 for the case where $x\in M^{sreg}(n)$ is
especially nice. \vs {\bf Theorem 0.8.} {\it Let $x\in M^{sreg}(n)$. Then the morphism (0.5) is an algebraic
isomorphism onto its image, the maximal orbit $A\cdot x$. In particular $A\cdot x$ is a nonsingular variety
and as such $$A\cdot x \cong G_{x,1}\times
\cdots \times G_{x,n-1}\eqno (0.7)$$}\vs See Theorem 3.14. \vs 0.3. Let $x\in M(n)$. Motivated by the Jacobi
matrices which arise in the theory of orthogonal polynomials on $\Bbb R$, we will say that
$x$ satisfies the eigenvalue disjointness condition if, for any $m\in I_n$, the eigenvalues of $x_m$ have
multiplicity one (so that $x_m$ is regular semisimple in $M(m)$) and, as a set, $E_{x}(m)\cap E_{x}(m+1)=
\emptyset$ for any
$m\in I_{n-1}$. Let $M_{\Omega}(n) $ be the dense Zariski open set of such $x\in M(n)$. One readily has
that
$M_{\Omega}(n) = \Phi_n^{-1}(\Omega(n))$ where $\Omega(n)$ is a dense Zariski open set in $\Bbb C^{d(n)}$.
\vs {\bf Theorem 0.9.} {\it One has $M_{\Omega}(n)\s M^{sreg}(n)$. In fact if $c\in \Omega(n)$ then the
entire fiber $M_c(n)$ is a single maximal $A$-orbit. Moreover if $c\in \Omega(n)$ and $x\in M_c(n)$ then
$G_{x,m}$ is a maximal (complex) torus in $Gl(m)$, for any $m\in I_{n-1}$, so that $M_c(n) = A\cdot x$ is a
closed nonsingular subvariety of $M(n)$ and as such $$M_c(n)\cong (\Bbb C^{\times})^{d(n-1)}\eqno
(0.8)$$}\vs See (2.55) and Theorem 3.23. Let $\u$ be the $d(n-1)$-dimensional nilpotent Lie algebra of
all strictly upper triangular matrices. One has a natural projection $$M(n)\to \u,\qquad x\mapsto
x_{\u}\eqno (0.9)$$ where $x_{\u}$ is such that $x-x_{\u}$ is lower triangular. Another ``snapshot" picture
of $M_c(n)$, when $c\in \Omega(n)$, is given in \vs {\bf Theorem 0.10.} {\it Let $c\in \Omega(n)$ and let
$\u_c$ be the image of $M_c(n)$ by the projection (0.9). Then $\u_c$ is a dense Zariski open subset of $\u$
and the restriction $$M_c(n)\to\u_{c}\eqno (0.10)$$ of (0.9) to $M_c(n)$ is an algebraic isomorphism.}\vs
See Theorem 3.26. 

Let $c\in \Omega(n)$. Then it follows from Theorem 0.9 that the subgroup $D_c = \{a\in A\mid a\cdot x =
x,\,\,\forall x\in M_c(n)\}$ is closed and discrete in $A$. Let
$A_c = A/D_c$ so that the action of
$A$ on $M_c(n)$ descends to an action of $A_c$. Furthermore the latter action is simple and transitive so
that 
$A_c$ has the structure of an algebraic group and as such $$A_c\cong (\Bbb C^{\times})^{d(n-1)}\eqno (0.11)$$
In particular if $F_c = \{a\in A_c\mid a^2 = 1\}$ then $F_c$ is a finite group of order $2^{d(n-1)}$. 

It is suggestive from the example of symmetric Jacobi matrices that if $x$ is a symmetric
matrix then
$x$ is close to being determined by knowing the spectrum of $x_m$ for all $m\in I_n$. Let $c\in
\Omega(n)$ and let $M_c^{(sym)}(n)$ be the set of all symmetric matrices in $M_c(n)$. The following
theorem appears in the paper as Theorem 3.32. \vs {\bf Theorem 0.11.} {\it Let $c\in \Omega(n)$ (see (2.53)).
Then
$M_c^{(sym)}(n)$ is a finite set of cardinality
$2^{d(n-1)}$. In fact $M_c^{(sym)}(n)$ is an orbit of
$F_c$.}\vs If $c\in \Bbb C^{d(n)}$ put $E_c(m) = E_x(m)$, for $m\in I_n$, and any (and hence all) $x\in
M_c(n)$. Also put $\mu_{k\,m}(c) = \mu_{k\,m}(x)$ for any $k\in I_m$. We will say that
$c\in
\Bbb C^{d(n)}$ satisfies the eigenvalue interlacing condition if one has $\mu_{k\,m}(c) \in \Bbb R$ for all
$m\in I_n,\,k\in I_m$ and, writing $\mu_{k\,m} = \mu_{k\,m}(c)$, 
$$\mu_{1\,m+1}<\mu_{1\,m}<
\mu_{2\,m+1}<\ldots < \mu_{m\,m+1}<\mu_{m\,m} <
\mu_{m+1\,m+1}\eqno (0.12)$$ for any $m\in I_{n-1}$. If
$c$ satisfies the eigenvalue interlacing condition then obviously
$c\in
\Omega(n)$. The following result is established in the paper as Theorem 3.34. \vs {\bf Theorem 0.12.}
 {\it Let $c\in \Omega(n)$. Then the following conditions are equivalent,
$$\eqalign{&(1)\,\,\,\hbox{$c$ satisfies the eigenvalue interlacing condition},\cr &(2)
\,\,\,\hbox{there exists a real
symmetric matrix in
$M_c^{(sym)}(n)$},\cr  &(3)\,\,\,\hbox{all $2^{d(n-1)}$ matrices in $M_c^{(sym)}(n)$ are real
symmetric}.\cr}$$}\vs {\bf Example 0.13.}  Consider the case where $n=3$ so that there are
8 symmetric matrices in $M_c(3)$ for any $c\in \Omega(n)$. Let $c$ be defined so that $E_c(1) = \{0\},\,E_c(2)
= \{1,-1\},\,E_c(3) = \{\sqrt 2,0,-\sqrt 2\}$ so that $c$ is eigenvalue interlacing.
Then 2 of the 8 real symmetric matrices in $M_c(3)$ are 
$$x=\left(\matrix{0&1&0\cr 1&0&1\cr 0&1&0\cr}\right),\qquad
y=\left(\matrix{0&1&1\cr 1&0&0\cr 1&0&0\cr}\right)$$ noting that $x$, but not $y$, is Jacobi. 
The remaining 6 are
obtained by sign changes in $x$ and $y$. \vs Let $\phi_k(t),\,k\in \Bbb Z_+$, be an orthonormal sequence of
polynomials in
$L_2(\Bbb R,\nu)$ (for a suitable measure $\nu$ on $\Bbb R$) obtained by applying the Gram--Schmidt process
to the monomial functions
$\{t^k\},m\in
\Bbb Z_+$. The $\phi_k$ are uniquely determined up to sign. Let $W_k$ be the span of $1,\ldots,t^{k-1}$ and
let
$\Lambda_k$ be the orthogonal projection of
$L_2(\Bbb R,\nu)$ on $W_k$ so that $\{\phi_{j-1}(t)\},\,j\in I_{k}$, is an orthonormal basis of $W_k$. Let 
$\tilde {t}$ be the operator on $L_2(\Bbb R,\nu)$ of multiplication by $t$. \vs {\bf Theorem 0.14.} {\it Let
$x\in M(n)$ be the matrix of $\Lambda_n\,\tilde {t}$ with respect to the basis $\phi_{k},\,k = 0,\ldots,
n-1$ and let $c = \Phi_n(x)$ (see (0.1)). Then $c$ satisfies the eigenvalue interlacing condition and $x\in
M_c^{(sym)}$. Moreover $x$ is Jacobi (see \S 2.5) and in the $2^{d(n-1)}$-element set $M_c^{(sym)}$ there are
precisely $2^{n-1}$ Jacobi matrices. The latter represents $\Lambda_n\,\tilde {t}$ when the basis $\phi_k$
is replaced, using sign changes, by $2^{n-1}$ different choices of the orthonormal polynomials. 

Finally the
characteristic polynomial of $x_m$, for $m\in I_n$, is the monic polynomial corresponding to $\phi_{m}$.
In particular the numbers in $E_x(m)$ are the zeros of the orthogonal polynomial $\phi_{m}$.}\vs See Theorems
2.20, 2.21 and (2.67) where the $a_k>0$.\vs 0.4. For any $x\in M(n)$ let $O_x$ be the $Gl(n)$ adjoint
orbit generated by $x$. The symplectic leaves of the Poisson structure
on
$M(n)$ are, of course, all the adjoint orbits and any adjoint orbit stable under the action of $A$. For any
$x\in M(n)$ let
$O^{sreg}_x = O_x\cap M^{sreg}(n)$ so that $O^{sreg}_x$ is an $A$-stable Zariski open set in $O_x$. \vs
{\bf Theorem 0.15.} {\it Let $x\in M(n)$. Then $O^{sreg}_x$ is not empty if and only if $x$ is regular. In
such a case $O^{sreg}_x$ is a $2\,d(n-1)$-dimensional symplectic (in the complex sense) open and dense
submanifold of $O_x$. Furthermore the orbits of $A$ in $O^{sreg}_x$ 
(necessarily of dimension $d(n-1)$) are the leaves of a polarization of $O^{sreg}_x$. }\vs See Theorem 3.36.
\vs Part II of this 2 part paper, to appear elsewhere, is 
devoted to solving a ``classical" analogue of the Gelfand--Kirillov conjecture. The main result there will
use results in the present paper. The eigenvalue functions $\mu_{k\,m}(x)$, for $m\in I_n$ and $k\in I_m$,
can only be defined locally on
$M(n)$. However they can be defined globally on a suitable covering of $M_{\Omega}(n)$ and they Poisson
commute on this covering. If $\b_{\Omega} = \b_e\cap M_{\Omega}(n)$ then $M_{\Omega}(n) = A\cdot 
\b_{\Omega}$. 
A second set of
$d(n-1)$ Poisson commuting functions on $M_{\Omega}(n)$ (and hence on the covering) is then defined using a
natural coordinate system on $A$. Taking certain quotients it is shown that these two sets of Poisson
commuting elements generate two sets of functions which satisfy the desired commutation (Heisenberg)
relations with each other. The generated functions are in an algebraic extension of the function field of
$M(n)$.\vs 0.5. We wish to thank Pavel Etingof for fruitful conversations.\vs

\rm

	\centerline{\bf 1. Preliminaries}\vskip 1pc 
       
1.1. For any positive integer $n$ let $M(n)$ be the algebra (both
Lie and associative) of all $n\times n$ complex matrices. With
respect to its Lie algebra structure let $U(n)$ be the universal
enveloping algebra of $M(n)$. Let $S(n) = \oplus_{k=0}^{\infty} S^k(n)$
be the (graded) symmetric algebra over $M(n)$. Then $S(n) = Gr\,U(n)$ with
respect to the usual filtration in $U(n)$. Commutation in $U(n)$ induces
the structure of a Poisson algebra $[u,v]$ on $S(n)$ where
$$[S^j(n),S^k(n)]\s S^{j+k-1}(n)\eqno (1.1)$$ If one identifies $M(n)$ with
$S^1(n)$ then the Poisson bracket on $M(n)$ induced by (1.1) is of
course the same as the given Lie algebra bracket. 

Let $P(n)= \oplus_{k=0}^{\infty} P^k(n)$ be the graded algebra of
polynomial functions on $M(n)$. We identify $P^1(n)$ with the dual space
$M(n)^*$ to $M(n)$. For any positive integer $k$ let $I_k = \{1,\ldots,k\}$. For $i,j\in I_n$ we 
will, throughout the paper, let $\alpha_{ij}\in P^1(n)$ be the linear function defined
so that, for $x\in M(n)$, $$\alpha_{ij}(x)\,\,\hbox{is the $ij^{th}$ entry of
$x$}\eqno (1.2)$$ 

Of course $M(n)$ is the Lie algebra of the general linear group $Gl(n)$. The
adjoint action of $Gl(n)$ on $M(n)$ induces a $GL(n)$-module structure on the
algebras $U(n),\,S(n)$ and $P(n)$. It will be convenient in this paper to transfer the
coadjoint orbit theory of $GL(n)$ to adjoint orbits and to transfer the Poisson
algebra structure on $S(n)$ to $P(n)$. Let $\gamma:M(n)\to P^1(n)$ be the linear
isomorphism defined by $$\langle\gamma(x),y\rangle = -B(x,y)\eqno (1.3)$$ where $B$ is the bilinear $B(x,y) =
tr\,xy$ on $M(n)$. Then 
$\gamma$ extends to an algebra $Gl(n)$-isomorphism
$$\gamma:S(n)\to P(n)\eqno (1.4)$$ Consequently $P(n)$ becomes a $Gl(n)$-Poisson
algebra where $[\gamma(u),\gamma(v)]= \gamma([u,v])$ for $u,v\in S(n)$. If $e_{ij}\in
M(n)$ is the $ij^{th}$ matrix unit of $M(n)$ one notes that $$\alpha_{ij} =
\gamma(-e_{ji})\eqno (1.5)$$ Since negative transpose defines a Lie algebra
automorphism on $M(n)$ the bracket structure for the linear functionals
 $\alpha_{ij}$ is the same as for the matrix units. That is,
$$\eqalign{[\alpha_{ij},\alpha_{st}] &= 0\,\,\hbox{if}\,\,i\neq t
\,\,\hbox{and}\,\,j\neq s \cr [\alpha_{ij},\alpha_{jk}] &=
\alpha_{ik}\,\,\hbox{if}\,\,i\neq k\cr [\alpha_{ij},\alpha_{ji}]&=
\alpha_{ii}-\alpha_{jj}\cr}\eqno (1.6)$$ 

1.2.  Let $T(M(n))$ and $T^*(M(n))$ be,
respectively, the holomorphic tangent and cotangent bundles of $M(n)$. All vector
fields and 1-forms considered here will be, respectively, holomorphic sections of
$T(M(n))$ and $T^*(M(n))$ defined on some open subset (most often on $M(n)$ itself) of
$M(n)$. The Poisson algebra structure on $P(n)$ defines a Poisson structure, in the
holomorphic sense, on
$M(n)$. See Chapter 1 in  [CG]. If $\varphi\in P(n)$ the map $\psi\mapsto
[\varphi,\psi]$ is a derivation of
$P(n)$ and hence there exists a vector field $\xi_{\varphi}$ on $M(n)$ such that
$$\xi_{\varphi}\psi = [\varphi,\psi]$$ Since $\xi_{\varphi}\psi = -\xi_{\psi}\varphi$ it is clear that
$(\xi_{\varphi})_x$, for any $x\in M(n)$, depends only on $(d\varphi)_x$ so that in fact one
has a bundle map $$\Gamma:T^*(M(n)) \to T(M(n))\eqno (1.7)$$ and $(\xi_{\varphi})_x
=\Gamma((d\varphi)_x)$. 

If $W\s M(n)$ is any open set let
${\cal H}(W)$ be the algebra of holomorphic functions on $W$ and let $Vec(W)$ be the
Lie algebra of (holomorphic) vector fields on $W$. For any $y\in M(n)$ let $\partial^y\in Vec(M(n))$ be
the translation invariant vector field on $M(n)$ defined so that if $\psi\in
{\cal H}(M(n))$ and 
$x\in M(n)$ then
$$\eqalign{(\partial^y\,\psi)(x) &= (\partial^y)_x \psi\cr &= d/dt (\psi(x
+ty))_{t=0}\cr}$$
Let $$M(n)\to Vec(M(n)),\qquad
y\mapsto \eta^y$$ be the corresponding Lie algebra homomorphism corresponding to the
adjoint action of $Gl(n)$. If
$O$ is an adjoint orbit then obviously
$\eta^y|O$ is tangent to $O$. Explicitly it is immediate that if $x,y\in M(n)$ then
$$(\eta^y)_x= (\partial^{-[y,x]})_x\eqno (1.8)$$ so that if $O_x$ is the adjoint
orbit containing $x$ then $$T_x(O_x) = \{(\partial^{-[y,x]})_x\mid y\in M(n)\}\eqno
(1.9)$$ Using $\gamma$ we may carry over the symplectic structure (in the
holomorphic sense, see [GC], Chapter 1) on coadjoint orbits to adjoint orbits. If
$x\in M(n)$ let
$\omega_x$ be the value of the symplectic form $\omega_{O}$ on $O= O_x$ at $x$. Then
for
$y,z\in M(n)$ (taking into account the minus sign in the definition of $\gamma$, and
(5.3.1) and (5.3.3) in [K]) one has
$$
\omega_x(\eta^y,\eta^z) = B(x,[y,z])\eqno (1.10)$$ in the notation (1.3). \vskip .5pc
{\bf Lemma 1.1.} {\it One has $$\xi_{\gamma(y)} = \eta^y\eqno (1.11)$$ for any $y\in
M(n)$.} \vskip 1pc {\bf Proof.} For any $v\in P^1(n)$ and $x\in M(n)$ one has
$$(\eta^y \,v)(x)  = -\langle v,[y,x]\rangle \eqno (1.12)$$ by (1.8). But then if
$w\in M(n)$ one has
$$(\eta^y \,\gamma(w))(x) = B(w,[y,x])\eqno (1.13)$$ On the other hand
$$\eqalign{(\xi_{\gamma(y)}\gamma(w))(x)&= [\gamma(y),\gamma(w)](x)\cr
&=\gamma ([y,w])(x)\cr &= B([w,y],x)\cr}$$ But then $(\eta^y \,\gamma(w))(x)=
(\xi_{\gamma(y)}\gamma(w))(x)$. This proves (1.11). QED\vs We may enlarge the family of
functions $\varphi$ for which the vector field $\xi_{\varphi}$ is defined. Let $W\s M(n)$
 be any open set. For any $\varphi\in {\cal H}$ let $\xi_{\varphi}$ be the vector field on
$W$ defined so that if $x\in W$ then $$(\xi_{\varphi})_x = \Gamma ((d\varphi)_x)\eqno 
(1.14)$$

Let
$y_j,\,j=1,\ldots,n^2$, be a basis of $M(n)$ and let $v_j\in P^1(M(n))$ be the dual
basis so that, in the notation of (1.14), $$d\varphi =
\sum_{j}\partial^{y_j}(f)\,dv_j\eqno (1.15)$$ on $W$. But this implies $$\xi_{\varphi} =
\sum_{j}\partial^{y_j}(f)\,\xi_{v_j}\eqno (1.16)$$ on $W$.\vs {\bf Remark 1.2.} Using (1.8),(1.9) and
(1.11) and the notation of (1.14), note that for any
$x\in W$ the tangent vectors
$(\xi_{v_j})_x,\,j=1,\ldots,n^2,$ span
$T_x(O_x)$. \vs The following proposition
essentially recovers the fact that the adjoint orbits are the symplectic leaves of the Poisson structure
on $M(n)$. It also removes a possible ambiguity in the definition of the vector fields $\xi_{\varphi}$.
If $x\in M(n)$ and $\psi$ is a holomorphic function defined in a neighborhood of $x$ in $O_x$ then the
Hamiltonian vector field $\xi_{\psi}$ defined in this neighborhood is such that
$$\omega_x((\xi_{\psi})_x,u) = u\,\psi\eqno (1.17)$$ for any $u\in T_x(O_x)$ (see (4.1.3) in [K]). \vs 
{\bf Proposition 1.3.} {\it Let $W\s M(n)$ be an open set and let $x\in W$. Let $\varphi\in {\cal H}(W)$.
Then $$(\xi_{\varphi})_x \in T_x(O_x)\eqno (1.18)$$ Furthermore if $V = W\cap O_x$ then
$$\xi_{\varphi}|V = \xi_{\varphi|V}\eqno (1.19)$$}\vs {\bf Proof.} One has (1.18) by Remark 1.2 and
(1.16). Let $\psi = \varphi|V$. Since $\omega_x$ is nonsingular, to prove (1.19) it suffices, by Remark
1.2 and (1.18), to prove $$\omega_x((\xi_{\varphi})_x,(\xi_{v})_x) =
\omega_x((\xi_{\psi})_x,(\xi_{v})_x)\eqno (1.20)$$ for any $v\in P^1(n)$. But now it follows from
(1.15), (1.16) and (1.17) that to prove (1.20) it suffices to assume $W=M(n)$ and to take $\varphi\in
P^1(n)$. Write $v = \gamma(y)$ and $\varphi=\gamma(w)$ for $y,w\in M(n)$. Then the right side of (1.20)
is $(\xi_{\gamma(y)}\gamma(w))(x)$ by (1.17). But the end of the proof of Lemma 1 implies that this equals
$tr\,[w,y]\,\,x$. On the other hand the left-hand side of (1.20) equals $tr\,x\,[w,y]$ by (1.10). This
proves (1.19). QED\vs  Classical properties (see e.g., (4.1.4) in [K] for the real case) of Poisson bracket of
functions and  Hamiltonian vector fields on symplectic manifolds remain valid for the Poisson manifold $M(n)$ .
\vs {\bf Proposition 1.4. } {\it Let
$W\s M(n)$ be any open set. Then ${\cal H}(W)$ is a Lie algebra under Poisson bracket and the map $${\cal
H}(W) \to Vec(W),\qquad \varphi\mapsto \xi_{\varphi}\eqno (1.21)$$ is a Lie algebra homomorphism.}\vs
{\bf Proof.} This is immediate from (1.19) since any adjoint orbit is a symplectic manifold. QED\vskip
1pc

{\centerline{\bf 2. Commuting vector fields arising from Gelfand--Zeitlin theory}\vskip 1pc 2.1. Let
$m\in I_n$ and let $L_{m}$ be the set of all $n^2-m^2$ pairs $(i,j)$ where $\{i,j\}\s I_n$ but $\{i,j\}\not\subset I_m$.
We identify (``upper left hand corner") 
$$M(m)=\{x\in M(n)\mid \alpha_{ij}(x) = 0,\,\,\forall (i,j)\in L_m\}$$ We will regard $M(m)$ as a Lie
subalgebra of $M(n)$. The corresponding subgroup $\{g\in Gl(n)\mid \alpha_{ij}(g)=\delta_{ij},\,\forall
(i,j)\in L_m\}$ naturally identifies with $Gl(m)$. One has the direct sum $$M(n) = M(m)\oplus M(m)^{\perp}\eqno
(2.1)$$ where $M(m)^{\perp}$ is the $B$-orthocomplement of $M(m)$ in $M(n)$. For any
$x\in M(n)$ let 
$x_m\in M(m)$ be the component of $x$ in $M(m)$ relative to (2.1). That is, $x_m\in M(m)$ is defined so that
$\alpha_{ij}(x_m) = \alpha_{ij}(x)$ for all
$(i,j)\in I_{m}$. We will also refer to $x_m$ as the $m\times m$ cutoff of $x$. The surjective map $M(n)\to
M(m)$ defined by (2.1) induces an injection $P^1(m)\to P^1(n)$. The latter extends to an injective homomorphism
$$P(m)\to P(n)\eqno (2.2)$$ It follows immediately from (1.6) that (2.2) is an injective homomorphism of Poisson
algebras and hence $M(m)$ embeds as a Poisson submanifold of $M(n)$. Henceforth we will identify $P(m)$ with
its image in $P(n)$ under (2.2).

Let $m\in I_n$ and let $Id_{m}$ be the identity element of the associative algebra $M(m)$. For any $k\in I_m$
let $f_{k,m}\in P(n)$ be the coefficient of the cutoff characteristic polynomial defined by $$det(\lambda\,\,
Id_m-x_m) = \sum_{k=1}^m\,(-1)^{m-k +1} f_{k,m}(x) \lambda^{k-1}\eqno (2.3)$$ One
readily notes that if $\{\mu_j\},\,j=1,\ldots,m$, are the eigenvalues (always with multiplicity as roots of the
characteristic polynomial) of
$x_m$, then
$$f_{k,m}(x)\,\,\hbox{is the elementary symmetric polynomial of degree $m-k+1$ in the $\mu_j$}\eqno (2.4)$$\vs For
any nonnegative integer
$k$ let
$d(k) = k(k+1)/2$. Note that $$d(n-1) + d(n) = n^2$$

Since $k,m$ are arbitrary subject to the condition $1\leq k\leq m\leq n$ the number of polynomials $f_{k,m}$
defined by (2.3) is 
$d(n)$. It will be convenient to simply order these polynomials, defining $p_j$, for $j\in
I_{d(n)}$, by putting
$$p_{d(m-1) + k} = f_{k,m}\eqno (2.5)$$ One readily deduces the following proposition from the commutativity of
the Gelfand--Zeitlin ring in the enveloping algebra $U(M(n))$. The proof given here is self-contained. \vs {\bf Proposition
2.1}. {\it One has the Poisson commutativity
$$[p_i,p_j] = 0\eqno (2.6)$$ for any $\{i,j\}\s I_{d(n)}$.} \vs {\bf Proof.} One of course knows that the
ring of invariants $P(n)^{Gl(n)}$ is the polynomial ring in $p_{d(n-1) +k},\,\,k\in I_n$. Hence for
these values of $k$ it follows that $p_{d(n-1) +k}$ is constant on any adjoint orbit $O_x$. This implies that
$$\xi_{p_{d(n-1) +k}} = 0\eqno (2.7)$$ by (1.19). Thus one has (2.6) if $j = d(n-1) +k$ and $k\in I_n$. But
then, replacing $n$ by any $m\leq n$, (2.6) follows for all $\{i,j\} \s I_{d(n)}$. QED\vs 2.2. Let $\Phi_n$ be
the regular algebraic map
$$\Phi_n:M(n)\to \Bbb C^{d(n)},\qquad \Phi_n(x) = (p_1(x),\ldots, p_{d(n)}(x))\eqno (2.8)$$ and for any
$c=(c_1,\ldots,c_{d(n)})\in \Bbb C^{d(n)}$ let $M_c(n) = \Phi_n^{-1}(c)$ so that $M_c(n)$ is a typical fiber of
$\Phi_n$ and hence $$M(n)= \sqcup_{c\in \Bbb C^{d(n)}}\,M_c(n)\eqno (2.9)$$ From the definition we note that
if
$x,y\in M(n)$ then $x,y$ lie in the same fiber of $\Phi_n$ if and only if the characteristic polynomial
of the cutoffs $x_m$ and $y_m$ are equal for all $m\in I_n$. Expressed otherwise, introduce a lexicographical order in $\Bbb
C$ so that if $z_1,z_2\in \Bbb C$ then $z_1< z_2$ if $\Re\, z_1< \Re\,z_2$ and $\Im\,z_1<\Im\,z_2$ in case $\Re\, z_1=
\Re\,z_2$. Let
$c\in
\Bbb C^{d(n)}$ and let, for $m\in I_n$,  $$ E_c(m) = \{\mu_{1\,m}(c),\ldots,\mu_{m\,m}(c)\}\eqno (2.10)$$ be the roots, in
increasing order, of the polynomial $$\lambda^m + \sum_{k=1}^m\,(-1)^{m-k +1} c_{d(m-1) +k}
\lambda^{k-1}\eqno (2.11)$$ Then $x\in M_c(n)$ if and only if $E_c(m)$ is an $m$-tuple of the eigenvalues of
$x_m$, for all $m\in I_n$, where the multiplicity is as a root of the characteristic polynomial of $x_m$.\vs {\bf Remark
2.2.} Of course given any set of
$n$ sequences, $ E(m) = \{\mu_{1\,m},\ldots,\mu_{m\,m}\},\, m\in I_n$, there exists a unique $c\in \Bbb C^{d(n)}$ such
that
$E(m) = E_c(m)$, up to an ordering, for all $m\in I_n$. \vs Let $e\in M(n)$ be the principal nilpotent
element $e = -\sum_{i=1}^{n-1}e_{i+1,i} $ and let $\b\s M(n)$ be the Borel subalgebra of all upper triangular
matrices. Thus $x\in e + \b$ if and only if $x$ is of the form  $$x =  \left(\matrix{a_{1\,1}&a_{1\,2}&\cdots
&a_{1\,n-1}&a_{1\,n}\cr -1&a_{2\,2}&\cdots &a_{2\,n-1}&a_{2\,n}\cr 0&-1&\cdots &a_{3\,n-1}&a_{3\,n}\cr
\vdots &\vdots &\ddots &\vdots &\vdots\cr 0&0&\cdots &-1&a_{n\,n}\cr}\right)\eqno (2.12)$$ where $a_{ij}\in \Bbb
C$ is arbitrary. That is $\alpha_{ij}(x) = a_{ij}$ for $i\leq j$. If one considers the $n$-dimensional
subvariety $\ss = \{x\in e+\b\mid a_{ij} = 0\,\,\hbox{for}\,\,1\leq i\leq j\leq n-1\}$ then from well--known
facts about companion matrices, for all $k\in I_n$ and $x\in \ss$,  $$\alpha_{k\,n}(x) =
p_{d(n-1) + k}(x)\eqno (2.13)$$ In particular the restriction of the $Gl(n)$-invariants, $p_{d(n-1) +
k}$ to $\ss$ define a coordinate system on $\ss$ and the map $$\ss\to \Bbb C^n,\qquad x\mapsto 
(p_{d(n-1) +1}(x),\ldots,p_{d(n-1) +n}(x))\eqno(2.14)$$ is an algebraic isomorphism. We now generalize this
statement to the entire affine variety $e+\b$. \vs {\bf Theorem 2.3.} {\it Let $q_i = p_i| e + \b$ 
for $i\in I_{d(n)}$. Then the $q_i$ are a coordinate system on $e+\b$. In particular for any $x\in e+\b$ the
differentials $$(dq_i)_x,\,i\in I_{d(n)},\,\,\hbox{are linearly independent}\eqno(2.15)$$ and hence
the differentials $$(dp_i)_x,\,i\in I_{d(n)},\,\,\hbox{are linearly independent}\eqno (2.16)$$ Furthermore the
map $$e + \b \to \Bbb C^{d(n)},\qquad x\mapsto (p_1(x),\ldots,p_{d(n)}(x))\eqno (2.17)$$ is an algebraic 
isomorphism.}\vs {\bf Proof.} We will prove Theorem 2.3 by induction on $n$. If $n=1$ the proof is immediate.
Assume that the theorem is true for $n-1$. Let $(e+\b)_{n-1}$ be equal to $e+\b$ when $n-1$ replaces $n$. Then
$x\mapsto x_{n-1}$ defines a surjection $$e+\b \to (e+\b)_{n-1}\eqno (2.18)$$ If $z\in (e+\b)_{n-1}$ and
$F(z)$ is the fiber of (2.18) over $z$ then $$F(z)\to \Bbb C^n,\qquad x \mapsto (\alpha_{1n}(x),\ldots,
\alpha_{nn}(x))$$ is an algebraic isomorphism. If $x\in F(z)$ then $x_j =z_j$ for $j=1,\ldots,n-1$.
 Now expanding the determinant $det(\lambda I_n - x)$ using the last column of $x\in F(z)$ one
has $$\eqalign{det(\lambda I_n - x)&= \lambda^n + \sum_{k=1}^{n}(-1)^{n-k+1} p_{d(n-1)+ k}(x)\lambda^{k-1}\cr 
& = (\lambda - a_{nn}(x))\,det\, (\lambda\,I_{n-1}-z) + \sum_{k=1}^{n
-1}(-1)^{n-k+1}\alpha_{kn}(x)\,det\,(\lambda\,I_{k-1}- z_{k-1})\cr
&= (\lambda - \alpha_{nn}(x))g_{n-1}(\lambda) +
\sum_{k=1}^{n-1}(-1)^{n-k+1}\alpha_{kn}(x)\,g_{k-1}(\lambda)\cr}\eqno (2.19)$$ where
$det\,(\lambda\,I_{0}- z_{0})$ is interpreted as 1 and $g_j(\lambda)$ is a monic polynomial of degree $j$ in 
$\lambda$
whose coefficients are fixed elements of $P(n-1)$ evaluated at $z$.  Comparing coefficients there exists a upper
trianglar
$n\times n$ matrix $b_{i,j}$ of fixed elements in $P(n-1)$ evaluated at $z$ such that
$$\eqalign{p_{d(n-1) + n}(x) &= \alpha_{n\,n}(x) + b_{n,n}\cr p_{d(n-1) + n-1}(x) &= \alpha_{n-1\,n}(x) +
b_{n-1,n-1} + b_{n-1,n}\alpha_{n\,n}(x)\cr p_{d(n-1) + n-2}(x) &= \alpha_{n-2\,n}(x) + b_{n-2,n-2} +
b_{n-2,n-1}\,\alpha_{n-1\,n}(x) + b_{n-2,n}\,\alpha_{n\,n}(x)\cr \vdots &=\vdots\cr
p_{d(n-1) +1}(x)&= \alpha_{1\,n}(x) + b_{1,1} + b_{1,2}\,\alpha_{2\,n}(x) +\cdots +
b_{1,n}\,\alpha_{n\,n}(x)\cr}\eqno (2.20)$$ By induction the $b_{i,j}$ are polynomials in $p_i(z)=p_i(x)$, for
$i\leq d(n-1)$. Given the triangular nature of (2.20) we can solve for the
$\alpha_{k\,n}(x)$ in terms of 
$p_j(x),\,j\in I_{d(n)}$ and hence we can write any $\alpha_{i\,j}|e+\b,\,i\leq j$ as a polynomial in $q_j$ for
$j\in I_{d(n)}$. But since $\alpha_{i\,j}|e+\b,\,i\leq j$ is a coordinate sysytem for $e+\b$, this proves (2.15)
(and hence also (2.16)) and the injectivity of (2.17). The surjectivity of (2.17) follows from induction and the
fact that the
$\alpha_{k\,n}(x)$ in (2.20) are arbitrary. QED \vs {\bf Remark 2.4.} Note that (2.16) in Theorem 2.3 implies that
the functions $p_j, \,j\in I_{d(n)}$, are algebraically independent in $P(n)$. Note also, upon conjugating $e+\b$
by an invertible diagonal matrix $g$, that Theorem 2.3 is still true if $e$ is replaced by any sum
$\sum_{i=1}^{n-1}z_i\,e_{i+1,i}$ where $z_i\in \Bbb C^{\times}$. This is clear, since $g$ is diagonal,
$$p_j(g\,xg^{-1}) = p_j(x)$$ for any $j\in I_{d(n)}$ and $x\in M(n)$. \vs Let $x\in M(n)$ and for $m\in I_n$ 
let $E_x(m) = \{\mu_{1,m}(x),\ldots,\mu_{m,m}(x)\}$ be the $m$-tuple of the eigenvalues of $x_m$, in increasing
(lexicographical) order, with multiplicity as roots of the characteristic polynomial. In the notation of (2.10), one
has
$$E_x(m) = E_c(m),\quad \forall m\in I_n,\,\,\iff x\in M_c(n) \eqno (2.21)$$ \vskip .5pc
Theorem 2.3 has a number of nice
consequences. For one thing it implies that
$e+\b$ is a cross-section of the surjection
$\Phi_n$. For another for $x\in e+\b$ the spectum of
$x_m$ over all $m\in I_n$ can be chosen arbitrarily and (independently) and in so doing uniquely determines $x$.
That is,
\vs {\bf Theorem 2.5.} {\it The restriction of the map $\Phi_n:M(n)\to \Bbb C^{d(n)}$ (see (2.8)) to $e+\b$
is an algebraic isomorphism $$e+\b\to \Bbb C^{d(n)}\eqno (2.22)$$ Furthermore given any set of $n$ sequences
$E(m)= \{\mu_{1\,m},\ldots,\mu_{m\,m}\},\,m\in I_n$, there exists a unique $x\in e+\b$ such that $E(m)= E_x(m)$, up to a
reordering, for all
$m\in I_n$. }\vs {\bf Proof.} One has only to note that if $c\in \Bbb C^{d(n)}$ then $$M_c(n)\cap (e + \b) = \{x\}\eqno
(2.23)$$ where $x\in e +\b$ is the unique point such that $p_i(x) = c_i$ for all $i\in I_{d(n)}$. The final statement follows
from Remark 2.2. QED\vs 2.3. For any $x\in M(n)$ let
$M(n)_x$ be the centralizer of $x$ in $M(n)$ and let $Z_{x,n}$ be the (necessarily commutative) associative
subalgebra generated by $x$ and $Id_n$. One knows $$\eqalign{&Z_{x,n} =
cent\,\,M(n)_x\cr dim\,&Z_x\leq n\leq
dim\,M(n)_x\cr}\eqno (2.24)$$ One says that $x$ is regular if $dim\,M(n)_x = n$ and as one knows $$x\,\,\hbox{is
regular}\,\,\iff\,\,dim\,Z_{x,n} = n\,\,\iff Z_{x,n} = M(n)_x\eqno (2.25)$$ If $x\in M(n)$ it follows from (2.24)
that
$dim\,O_x
\leq 2\,d(n)$ and $$x\,\,\hbox{is regular}\,\,\iff\,\,dim\,O_x = 2\,d(n)\eqno (2.26)$$ Let $T_x(O_x)^{\perp}$ be
the orthocomplement of $T_x(O_x)$ in $T^*_x(M(n)$. \vs {\bf Proposition 2.6.} {\it Let $x\in M(n)$. Then $x$ is
regular if and only if $\{(df_{k,n})_x\},\,k\in I_n,$ are linearly independent. In fact $x$ is regular if and
only if $\{(df_{k,n})_x\},\,k\in I_n,$ is a basis of $T_x(O_x)^{\perp}$.} \vs {\bf Proof.} The first
statement follows from Theorem 9, p. 382 in [K-1]. If $\{(df_{k,n})_x\},\,k\in I_n,$ is a basis of
$T_x(O_x)^{\perp}$ then necessarily $dim\,O_x = 2\,d(n)$ so that $x$ is regular by (2.26). Conversely
$(df_{k,n})_x\in T_x(O_x)^{\perp}$ for any $k\in I_n$ and any $x\in M(n)$ since $f_{k,n}$ is constant on $O_x$.
But then $\{(df_{k,n})_x\},\,k\in I_n,$ is a basis of
$T_x(O_x)^{\perp}$ by the first statement and (2.26). QED\vs We will say $x\in M(n)$ is strongly regular if
$(dp_i)_x$ for all $i\in I_{d(n)}$ are linearly independent. Let $M^{sreg}(n)$ be the set of all strongly
regular elements in $M(n)$. By (2.16) one has $$e+\b\s M^{sreg}(n)\eqno (2.27)$$ It is then clear that
$$M^{sreg}(n)\,\,\hbox{is a nonempty Zariski open subset of $M(n)$}\eqno (2.28)$$ \vs {\bf Theorem 2.7.} {\it
Let $x\in M(n)$ so that $(\xi_p)_x\in T_x(O_x)$ by (1.18) for any $p\in P(n)$. Then $x$ is strongly regular if and
only if the tangent vectors $(\xi_{p_i})_x\in T_x(O_x),\,i\in I_{d(n-1)},$ are linearly independent.} \vs {\bf
Proof.} Assume that $x$ is strongly regular. Then $x$ is regular by Proposition 2.6. But if $q_i = p_i\mid O_x$
then $(dq_i)_x$, for all $i\in I_{d(n-1)}$ is a linearly independent set since otherwise there is a nontrivial
linear combination of $(dp_i)_x$ for $i\in I_{d(n-1)}$ which lies in $T_x(O_x)^{\perp}$. But this and
Proposition 2.6 would contradict the strong regularity of $x$. But now $$\xi_{q_i} = \xi_{p_i}|O_x\eqno (2.29)
$$
by (1.19). But this proves the
tangent vectors $(\xi_{p_i})_x\in T_x(O_x),\,i\in I_{d(n-1)},$ are linearly independent since $O_x$ is
symplectic. Conversely assume that the
tangent vectors $(\xi_{p_i})_x\in T_x(O_x),\,i\in I_{d(n-1)},$ are linearly independent. Then recalling (1.7)
it follows that the covectors $(dp_i)_x$ for $i\in I_{d(n-1)}$ are linearly independent. But then, by (2.29),
the same statement is true if we replace the $p_i$ by the restrictions $q_i$ (using the notation of (2.29)). Let
$W_x\s T_x(O_x)$ be the span of $(\xi_{p_i})_x\in T_x(O_x),\,i\in I_{d(n-1)}$ so that $dim\,W_x = d(n-1)$. But
(2.6) and (1.19) imply that $\omega_x((\xi_{p_i})_x,(\xi_{p_j})_x)= 0$ for $i,j\in I_{d(n-1)}$. Thus $W_x$ is a
Lagrangian subspace of the symplectic vector space $T_x(O_x)$. Thus $dim\,O_x \geq 2\,d(n-1)$. Hence $x$ is
regular. Thus $(dp_{d(n-1) + k})_x$ for $k\in I_{n}$ are linearly independent by Proposition 2.6. But then the
full set $(dp_i)_x$, for $i\in I_{d(n)}$, is linearly independent by Proposition 2.6 and the fact that the $(dq_i)_x$,
for $i\in I_{d(n-1)}$, are linearly independent. QED\vs 2.4. For $m\in I_n$ one knows that $P(m)^{Gl(m)}$ is a
polynomial ring with the $m$ generators $p_{d(m-1) + k},\,k\in I_m$. Let $J(n)\s P(n)$ be the algebra generated by
$P(m)^{Gl(m)}$, over all $m\in I_n$. By Remark 2.4 $$J(n)\cong P(1)^{Gl(1)}\otimes \cdots \otimes
P(n)^{Gl(n)}\eqno (2.30)
$$
and
$J(n)$ is the polynomial ring with the $d(n)$ generators $p_i,\,i\in I_{d(n)}$. It then follows from (2.6) that
$$[p,q] = 0,\,\,\hbox{for all $p,q\in J(n)$}\eqno (2.31)$$ and hence for the corresponding vector fields
$$[\xi_p,\xi_q]= 0 ,\,\,\hbox{for all $p,q\in J(n)$}\eqno (2.32)$$ (see Proposition 1.4). Now for any $x\in M(n)$
one defines a subspace of $T_x(M(n))$ by putting $$V_x = \{(\xi_p)_x\mid \forall p\in J(n)\}$$ \vskip .5pc {\bf
Remark 2.8.} If $x\in M(n)$ note that it follows from Proposition 1.3, (1.19) and (2.32) that $V_x\s 
T_x(O_x)$ and $V_x$ is an isotropic subspace of $T_x(O_x)$ with respect to the symplectic form $\omega_x$.
In particular
$$dim\,V_x\leq {1\over 2}\,dim\,O_x\eqno (2.33)$$ and hence $$dim\,V_x \leq d(n-1)\eqno (2.34)$$ Moreover note
that if $J(n)$ is generated by $\{q_i\},\,i\in I_k$, for some integer $k$ then $V_x$ is spanned by
$(\xi_{q_i})_x,\,i\in I_k$ since, clearly, $(dp)_x$ is in the span of $(dq_i)_x,\,i\in I_k$ (see (1.7)) for any
$p\in J(n)$. Recalling (2.7) note then that Theorem 2.7 may be expressed by the statement $$\hbox{One has equality
in (2.34) $\iff\,\, x$ is strongly regular}\eqno (2.35)$$\vs {\bf Remark 2.9.} Even though the function 
$x\mapsto dim\,V_x$ on $M(n)$ is nonconstant by, abuse of terminology, we will refer to $x\mapsto V_x$ as a distribution on
$M(n)$ (one readily shows that the two other conditions in the definition of an involutory distribution are satisfied). An
analytic submanifold $X$ of $M(n)$ will be said to be a leaf of this distribution if $T_x(X) = V_x$ for all $x\in X$.\vs 
One of main results of the paper will be to show that that there exists a connected analytic group $A$, operating on $M(n)$,
all of whose orbits are leaves of the distribution $x\mapsto V_x$. In addition it will be shown that $\xi_p$, for any $p\in
J(n)$, integrates to a flow on $M(n)$. Furthermore the flow commutes with $A$ and stabilizes the orbits of $A$. In  order to
do this we will first determine
$V_x$ very explicitly. For any
$m\in I_n$ and
$x\in M(n)$ let $V_x(m) = \{(\xi_{p})_x\mid p\in P(m)^{Gl(m)}\}$. It is immediate (2.7) that $V_x(n)= 0$ and from Remark 2.8
one has (choosing $k= d(n)$ and
$q_i=p_i$) that $$V_x = \sum_{m\in I_{n-1}} V_x(m)\eqno (2.36)$$ For $m\in I_n$ let $M(m)^{\perp}$ be defined as
in (2.1). Clearly if $u\in M(m)$ then
$$u + M(m)^{\perp} = \{y\in M(n)\mid y_m = u\}\eqno (2.37)$$ To determine
$V_x(m)$ we replace the generators
$f_{k,m}$ (see (2.3)) of $P(m)^{Gl(m)}$ by the generators $f_{(k,m)}$ of $P(m)^{Gl(m)}$ (recalling the theory of
symmetric functions) where for $k\in I_m$ and $x\in M(n)$, $$f_{(k,m)}(x) = {1\over m+1-k}\,\,tr\,(x_{m})^{m+1-k}\eqno
(2.38)$$
\vskip .5pc {\bf Lemma 2.10.} {\it Let $x\in M(n)$ and let $m\in I_n,\,k\in I_m$. Then for any $v\in M(n)$ and
any $w\in x_m + M(m)^{\perp}$ one has $$(\partial^v)_w f_{(k,m)} = tr\,(x_m)^{m-k}\,v\eqno (2.39)$$}\vs {\bf
Proof.} By definition $$\eqalign{(\partial^v)_w f_{(k,m)} &= {1\over m+1-k} \,\,d/dt(tr\,((w + tv)_m)^{m+1-k})_{t=0}\cr
& = {1\over m+1-k}\,\, d/dt(tr\,((x_m + tv_m)^{m+1-k})_{t=0}\cr
&= tr\,(x_m)^{m-k}\,v_m\cr
&= tr\,(x_m)^{m-k}v\cr}$$ since $(x_m)^{m-k}\in M(m)$ and hence $(x_m)^{m-k}$ is $B$-othogonal to the component
of
$v$ in $M(m)^{\perp}$ relative to the decomposition $M(n) = M(m) + M(m)^{\perp}$. QED\vs Let $x\in M(n)$ and let
$$T^*_x(M(n))\to P^1(n),\qquad \rho \mapsto \widetilde {\rho}\eqno (2.40)$$ be the isomorphism defined so that if $v\in
M(n)$ and $\rho \in T^*_x(M(n))$ then $\langle \rho,(\partial^v)_x\rangle = \langle\widetilde {\rho},v\rangle$. Clearly
$(d\widetilde {\rho})_x = \rho$ so that if $q\in P(n)$ and we put $\rho = (dq)_x$  then $$(\xi_{\tilde {\rho}})_x =
(\xi_q)_x\eqno (2.41)$$ by (1.7). \vs {\bf Proposition 2.11.} {\it Let $x\in M(n)$. Let $m\in I_n$ and $k\in I_m$. Let
$\gamma$ be defined as in (1.3). Then $$\gamma((x_m)^{m-k}) = - \widetilde {(df_{(k,m)})_x}\eqno (2.42)$$}\vs {\bf Proof.}
This is immediate from (2.39) since we can put $w=x$ in (2.39). QED \vs For notational convenience put $\xi_{(k,m)} =
\xi_{f_{(k,m)}}$. \vs {\bf Theorem 2.12.} {\it Let $m\in I_n$ and $k\in I_m$. Then for any $x\in M(n)$ one has
$$(\xi_{(k,m)})_x = (\partial^{[(x_m)^{m-k},x]})_x\eqno (2.43)$$}\vs {\bf Proof.} Put $y = -(x_m)^{m-k}$.
Then by (1.11) and (1.8) one has $$(\xi_{\gamma (y)})_x = (\partial^{-[y,x]})_x\eqno (2.44)$$ But then if $\rho = 
(df_{(k,m)})_x$ one has $$\gamma(y) = \widetilde {\rho}\eqno (2.45)$$ by (2.42). But 
$(\xi_{\tilde{\rho}})_x =(\xi_{(k,m)})_x $ by (2.41) where $q= f_{(k,m)}$. But then (2.43)
follows from (2.44) and (2.45). QED\vs For any $x\in M(n)$ and $m\in I_n$ let $Z_{x,m}$ be the associative
subalgebra generated by $Id_m$ and $x_m$ so that $Z_{x,m}$ is spanned by $(x_m)^{m-k},\,k\in I_m$. Here $(x_m)^0 =
Id_m$. Also let
$M(m)_{x_m}$ be the centalizer of
$x_m$ in
$M(m)$. One has that $dim\,Z_{x,m}\leq m$ and $x_m$ is regular in $M(m)$ if and only if $dim\,Z_{x,m} = m$ 
which is the case if and only if $$Z_{x,m} = M(m)_{x_m}\eqno (2.46)$$ Also put $$Z_x = \sum_{m\in I_n} Z_{x,m}\eqno
(2.47)$$ As a corollary of Theorem 2.12 one has
\vs {\bf Theorem 2.13}. {\it Let $x\in M(n)$. Then for any  $m\in I_m$, $$\eqalign{V_x(m) &=
\{(\partial^{[z,x]})_x\mid z\in Z_{x,m}\}\,\,\,\hbox{and}\cr V_x&= \{(\partial^{[z,x]})_x\mid
z\in Z_x\}\cr} \eqno (2.48)$$ Moreover $x$ is strongly regular if and only if (1) $x_m$ is regular in $M(m)$ for
all $m\in I_n$ and (2) the sum (2.47) is a direct sum. Equivalently $x$ is strongly regular if and only if the 
elements $(x_m)^{m-k}$, over all $m\in I_n$ and $k\in I_m$, are linearly independent in $M(n)$.} \vs {\bf Proof.} 
The equalities (2.48) follow from (2.36) and of course Theorem 2.12. See also the last part of Remark 2.8. The
statement about strong regularity follows from Proposition 2.11 and the isomorphism (2.41). QED\vs The following
result is a simpler criterion for strong regularity. \vs {\bf Theorem 2.14.} {\it
Let $x\in M(n)$. Then $x$ is strongly regular if and only if (a) $x_m$ is regular in $M(m)$ for all $m\in I_n$ and
(b) $$Z_{x,m}\cap Z_{x,m+1} = 0,\,\,\,\,\,\forall\,m\in I_{n-1}\eqno (2.49)$$} \vs {\bf Proof.} If $x$ is strongly
regular then (a) and (b) are satisfied by (1) and (2) in Theorem 2.13. Conversely assume (a) and (b) are satisfied.
We have only to show that (2) in Theorem 2.13 is satisfied. Assume (2) is not satisfied and so that there exists
$0\neq z_i\in Z_{x,m_i},\,i\in I_k$, for $k>1$, and $m_i$ a strictly increasing sequence with values in $I_n$, such
that $$\sum_{i=1}^k z_i = 0\eqno (2.50)$$ Let $m= m_1$ so that $m\in I_{n-1}$. But now if $i>1$ then note that
$$[z_i,x]_{m+1} =0\eqno (2.51)$$ Indeed $$M(n) = M(m_i) \oplus M(m_i)^{\perp}\eqno (2.52)$$ by (2.1). But $z_i\in
M(m_i)$ and clearly both components in (2.52) are stable under $ad\,M(m_i)$. However the component of $x$ in
$M(m_i)$ is
$x_{m_i}$. But $z_i\in Z_{x,m_i}$ so that $z$ commutes with $x_{m_i}$. Thus $[z_i,x]\in M_{m_i}^{\perp}$. But
$m_i\geq m+1$ and hence $M_{m_i}^{\perp}\s M_{m+1}^{\perp}$. This proves (2.51). But then (2.50) implies that
$[z_1,x]_{m+1} = 0$. But $ad\,M(m)$ clearly stabilizes both components of the decomposition $M(n) = M(m+1)\oplus
M(m+1)^{\perp}$. Since the component of $x$ in $M(m+1)$ is $x_{m+1}$ one has $[z_1,x_{m+1}] = 0$. Thus $z_1\in
M(m+1)_{x_{m+1}} = Z_{x,m+1}$. But then $z_1\in Z_{x,m}\cap Z_{x,m+1}$. This is a contradiction since $z_1\neq 0$.
QED\vs Recall (2.10). We will say that $c\in \Bbb C^{d(n)}$ satisfies
the eigenvalue disjointness condition if the (I) the numbers in $E_c(m)$  are distinct for all $m\in I_n$ and (II),
as a set, one
has
$E_c(m)\cap E_c(m+1) = \emptyset$ for all $m\in I_{n-1}$. Similiarly we will say $x$ satisfies the eigenvalue

disjointness condition if (I) and (II) are satisfied where $c$ is replaced by $x$. See (2.21). As explained later
conditions (I) and (II) are motivated by the theory of orthogonal polynomials. If $x$ satisfies the eigenvalue
disjointness condition, by abuse of notation, we will, on occasion, regard $E_x(m)$ as a set, rather than as a sequence. \vs
{\bf Remark 2.15}. Note that if
$x\in M(n)$ satisfies the eigenvalue disjointness condition then $x_m$ is a regular semisimple element of $M(m)$ for all
$m\in I_n$. \vs Let
$\Omega(n)$ be the set of all $c\in \Bbb C^{d(n)}$ which satisfy the eigenvalue disjointness condition and let
$M_{\Omega}(n)$ be the set of all
$x \in M(n)$ which satisfy this condition so that in the notation of (2.8), $$M_{\Omega}(n) =
\Phi_n^{-1}(\Omega(n))\eqno (2.53)$$\vs {\bf Remark 2.16}. We note here that $\Omega(n)$ and $M_{\Omega}(n)$ are
(obviously nonempty) Zariski open subsets of $\Bbb C^{d(n)}$ and $M(n)$, respectively. Indeed if
$c\in
\Bbb C^{d(n)}$ then condition I is just the condition that the discriminant of the polynomial (2.11) is
nonzero for all $m\in I_n$. Condition II is that
$$\prod_{i\in I_{m},j\in I_{m+1}} \mu_{i,m}(c)-\mu_{j,m+1}(c)\eqno (2.54)$$ not vanish for $m\in I_{n-1}$. But the
polynomial functions (2.54) on the roots $E_c(m),\,m\in I_n$, are clearly invariant under the product of the
root permutation groups $\Sigma(m),\,m\in I_n$. Consequently (2.54) can be expressed as polynomial functions on
$\Bbb C^{d(n)}$. Thus conditions I and II are the nonvanishing of certain polynomial functions on $\Bbb C^{d(n)}$.
A similiar statement clearly then applies to $M(n)$ thereby establishing the assertion of Remark 2.16.\vs 
Although with a different application in mind the constructions in the proof of the following theorem appear already
in
\S 4  of [GS]. \vs {\bf Theorem 2.17.} {\it If $x\in M(n)$ satisfies the eigenvalue disjointness condition then $x$
is strongly regular. That is, $$M_{\Omega}(n)\s M^{sreg}(n)\eqno (2.55)$$}\vs {\bf Proof.} Assume $x\in M(n)$ 
satisfies the eigenvalue disjointness condition. Let $m\in
I_{n-1}$. It suffices by (2.46) and Theorem 2.14 to prove that $$Z_{x,m}\cap M(m+1)_{x_{m+1}}= 0 \eqno (2.56)$$ But
now since
$x_m$ is regular semisimple it can be diagonalized. That is, there exists
$g\in Gl(m)$ such that $$g\, x_m\,g^{-1} = diag(\mu_{1\,m}(x),\ldots,\mu_{m\,m}(x))\eqno (2.57)$$ Then, writing
$\mu_{i\,m} = \mu_{i\,m}(x)$, $g\, x_{m+1}\,g^{-1}$ is of the form $$g\, x_{m+1}\,g^{-1} = 
\left(\matrix{\mu_{1\,m}&0&\cdots &0&a_{1\,m+1}\cr 0&\mu_{2\,m}&\cdots &0&a_{2\,m+1}\cr
\vdots &\vdots &\ddots &\vdots &\vdots\cr 0&0&\cdots &\mu_{m\,m}&a_{m\,m+1}\cr
b_{m+1\,1}&b_{m+1\,2}&\cdots &b_{m+1\,m}&d\cr}\right)\eqno (2.58)$$ But now 
$\{\mu_{1\,m+1}(x),\ldots,\mu_{m+1\,m+1}(x)\}$ are roots of the characteristic
polynomial $p(\lambda)$ of $x_{m+1}$. On the other hand, by (2.58) $$p(\lambda) = (\lambda-d)(\prod_{i=1}^m
(\lambda-\mu_{i\,m})  - \sum_{i=1}^m a_{i\,m+1}b_{m+1\,i}\prod_{j\in (I_m - \{i\})}(\lambda-\mu_{j\,m})\eqno
(2.59)$$ But $\mu_{i\,m}$ is a root of all the $m+1$ polynomial summands of (2.59) except the summand  
$$a_{i\,m+1}b_{m+1\,i}\prod_{j\in (I_m - \{i\})}(\lambda-\mu_{j\,m})$$ Since, as sets, $E_x(m)\cap E_x(m+1) = \emptyset$ one
must have that $$a_{i\,m+1}b_{m+1\,i}\neq 0,\,\,\forall i\in I_m\eqno (2.60)$$ But now the subalgebra of $M(m)$ generated by
$Id_{m}$ and $g\,x_m g^{-1}$ is the algebra $\d(m)$ of all diagonal matrices in $M(m)$. But now if $y=
g\,x_{m+1}\,g^{-1}$ and $M(m+1)_{y}$ is the centralizer of $y$ in $M(m)$ to prove (2.56) it
suffices, by conjugation, to prove  $$\d(m)\cap M(m+1)_y = 0\eqno (2.61)$$ But if $0\neq d \in \d(m)$ and $d =
diag(d_1,\,\ldots,d_m)$ then $d_j\neq 0$ for some $j\in I_m$. But $$\alpha_{j\,m+1}([d,y]) = d_j\,a_{j\,m+1}\eqno
(2.62)$$ Thus $[d,y] \neq 0$ by (2.60). This proves (2.61). QED\vs 2.5. If $x\in M(n)$ and $x$ is a real
symmetric matix then $x_m$, for any $m\in I_n$, is real symmetric so that the eigenvalues of $x_m$ are real. 
Hence, for the $m$-tuple, $E_x(m)$, one has 
$$\mu_{1\,m}\leq \cdots \leq \mu_{m\,m}\eqno (2.63)$$ where we have written $\mu_{j,m}$ for $\mu_{j,m}(x)$. 
Under certain conditions strong regularity implies the eigenvalue disjointness condition. A key argument in the proof
of Proposition 2.18 is well known.
\vs {\bf Proposition 2.18.} {\it Assume
$x\in M(n)$ is strongly regular and $x$ is a real symmetric matrix, so that in the notation of (2.63), the
inequality $\leq$ is replaced by the strict inequality $<$. Then
$x$ satisfies the eigenvalue disjointness condition. In fact if $m\in I_{n-1}$ then using the notation of (2.63) one
has the interlacing  condition $$\mu_{1\,m+1}<\mu_{1\,m}< \mu_{2\,m+1}<\ldots < \mu_{m\,m+1}<\mu_{m\,m} <
\mu_{m+1\,m+1}\eqno (2.64)$$}\vs {\bf Proof.} We use the notation in the proof of Theorem 2.17. But now 
instead of assuming that
$x$ satisfies the eigenvalue disjointness condition assume that $x$ is strongly regular so that one has (2.56) and
hence also (2.61). Since
$x$ is symmetric the element $g\in Gl(m)$ can be chosen to be orthogonal. Thus one has $a_{i\,m+1}= b_{m+1,i}$ for all
$i\in I_m$. But (2.61) implies that $[d,y] \neq 0$. On the other hand clearly $\alpha_{ik}([d,y]) = 0$ for all
$i,k\in I_n$ unless, as an ordered pair, $\{i,k\} = \{j,m+1\}$ or $\{i,k\} = \{m+1,j\}$. But  $\alpha_{j\,m+1}
([d,y]) = d_j\,a_{j\,m+1}$ and $\alpha_{m+1\,j} ([d,y]) = -d_j\,a_{j\,m+1}$. Thus $a_{j\,m+1}\neq 0$ for all $j\in
I_m$. But as a rational function of the parameter $\lambda$ of $\Bbb R$ one has $$ p(\lambda)/(\prod_{i=1}^m
(\lambda-\mu_{i\,m}) = \lambda - d - \prod_{j\in
(I_m - \{i\})}a_{j\,m+1}^2/(\lambda-\mu_{i\,m})\eqno (2.65)$$ Consider $\lambda$ in the open interval
$(\mu_{i\,m},\mu_{i+1\,m})$ for $i\in I_{m-1}$. As $\lambda$ approaches $\mu_{i\,m}$ the function (2.65) clearly
approaches $-\infty$ and (2.65) approaches $+\infty$ as $\lambda$ approaches $\mu_{i+1\,m}$. Thus there exists 
$\mu_{i+1\,m+1}\in (\mu_{i\,m},\mu_{i+1\,m})$. On the other hand if we consider $\lambda$ in the open interval
$(-\infty, \mu_{1\,m})$ one notes that (2.65) approaches $-\infty$ as $\lambda$ approaches $-\infty$ and (2.65)
approaches $+\infty$ as $\lambda$ approaches $\mu_{1\,m}$. Thus there exists $\mu_{1\,m+1}\in (-\infty, \mu_{1\,m}
)$. Similarly there exists $\mu_{m+1,m+1}\in (\mu_{m\,m},+\infty)$. This accounts for $m+1$ distinct roots of
$p(\lambda)$ not only establishing $E_x(m)\cap E_x(m+1) = \emptyset$ (so that $x$ satisfies the eigenvalue 
disjointness condition) but also (2.64). QED\vs In this paper $x\in M(n)$ is called a Jacobi matrix if
$\alpha_{ij}(x) = 0$ if $|i-j|>1$ and $\alpha_{ij}(x) \neq 0$ if $|i-j| = 1$. That is, $x$ is ``tridiagonal" 
except
that the main diagonal is arbitrary and the entries of $x$ are nonzero along the two diagonals adjacent to the
main diagonal. If these nonzero entries are all positive we will say that $x$ has positive adjacent diagonals.
Among other statements the following result says that if
$x$ is real symmetric and Jacobi the assumption of strong regularity in Proposition 2.18 is unnecessary since it is
automatically satisfied. In particular the conclusions of that proposition hold. 

It is clear that $M_c(n)$ is stable under conjugation by any diagonal matrix in $Gl(n)$ for any $c\in \Bbb
C^{d(n)}$.
\vs {\bf Theorem 2.19.} {\it Assume $x\in M(n)$ is a Jacobi matrix.  Then $x$ is strongly regular. Furthermore if
$y\in M(n)$ is also a Jacobi matrix then $x$ and $y$ are conjugate by a diagonal matrix in $Gl(n)$ if and only if 
$E_x(m) = E_y(m)$ for all $m\in I_n$ (i.e. if and only if $x$ and $y$ lie in the same fiber $M_c(n)$
of 
$\Phi_n$ -- see (2.21)). Moreover if $x$ real symmetric then
$x$ satisfies the eigenvalue disjointness condition and one has the eigenvalue interlacing $$\mu_{1\,m+1}<\mu_{1\,m}<
\mu_{2\,m+1}<\ldots <
\mu_{m\,m+1}<\mu_{m\,m} <
\mu_{m+1\,m+1}\eqno (2.66)$$ where we have written $\mu_{j\,k} = \mu_{j\,k}(x)$.}\vs {\bf Proof.} The first
statement follows from (2.16) and Remark 2.4. Given Jacobi matrices $x,y$ it follows from Remark 2.5 that there
exists diagonal matrices $g_x$ and $g_v$ in $Gl(n)$ such that $g_x\,x\,g_x^{-1}$ and $g_y\,y\,g_y^{-1}$ both lie
in $e + \b$. But then $g_x\,x\,g_x^{-1} = g_y\,y\,g_y^{-1}$ by Theorem 2.5 if $x$ and $y$ lie in the same fiber
$M_c(n)$. Thus the same fiber condition implies that
$x$ and
$y$ are conjugate by a diagonal matrix in $Gl(n)$. The converse direction is obvious. But now if $x$ is real
symmetric then the remaining statements follow from Proposition 2.18. QED\vs Theorem 2.19 has an immediate
application to the theory of orthogonal polynomials. Let $\rho(t)$ be a nonnegative integrable (with respect to
Lebesgue measure) function on a (finite or infinite) interval $(a,b)\s \Bbb R$. Assume that the measure $\nu =
\rho(t) dt$ on
$(a,b)$ has a positive integral and $\phi(t)$ is integrable with respect to $\nu$ for any polynomial function
$\phi(t)$. Let $W_k$ for $k\in \Bbb Z_+$ be the $k$--dimensional subspace of the Hilbert space $L_2(\Bbb R,\nu)$
spanned by $1,\ldots,t^{k-1}$. Let $\phi_k(t),\,k\in \Bbb Z_+$, be the orthonormal sequence in $L_2(\Bbb R,\nu)$
obtained by applying the Gram--Schmidt process to the functions $\{t^k\},k\in \Bbb Z_+$, so that, for any 
$k\in \Bbb Z_+$, $\{\phi_{j-1}(t)\},\,j\in I_{k}$, is an orthonormal basis of $W_k$. This basis 
may be normalized so that $a_k>0$ where $a_k\,t^k$ is the leading term of $\phi_k(t)$. Then one
knows that there is a 3-term formula (see e.g. (10), p. 157 in [J]) $$t\, \phi_k(t) = {a_k\over
a_{k+1}}\,\phi_{k+1}(t) + c_{kk}\,\,\phi_k(t) + {a_{k-1}\over a_{k}}\phi_{k-1}(t)\eqno (2.67)$$ where $c_{kk}\in
\Bbb R$. Let $\Lambda_k$ be the orthogonal projection of $L_2(\Bbb R,\nu)$ onto $W_k$. Let $\tilde {t}$ be the
multiplication operator on
$L_2(\Bbb R,\nu)$ by the function $t$. Then if $x\in M(n)$ is the matrix of $\Lambda_n\,\tilde {t}| W_n$
with respect to the basis $\phi_{m-1}(t),\,m\in  I_{n}$, of $W_n$, it follows from (2.67) that $x$ is the real
symmetric Jacobi matrix given by $$x =  \left(\matrix{c_{00}&a_{0}/ a_{1}&0&0&\cdots
&0\cr a_{0}/ a_{1}&c_{11}& a_{1}/ a_{2} &0&\cdots &0\cr 0 &a_{1}/ a_{2}
&c_{22}&a_{2}/ a_{3}&\cdots &0\cr
\vdots&\vdots &\vdots &\ddots &\vdots &\vdots\cr 0&0&0&\cdots &a_{n-1}/a_{n}&c_{nn}\cr}\right)\eqno (2.68)$$ By 
Theorem 2.19 $x$ is strongly regular and satisfies the eigenvalue disjointness condition. In addition one has the
eigenvalue interlacing
$$\mu_{1\,m+1}<\mu_{1\,m}<
\mu_{2\,m+1}<\ldots <
\mu_{m\,m+1}<\mu_{m\,m} <
\mu_{m+1\,m+1}\eqno (2.69)$$ where we have put $\mu_{j\,k} = \mu_{j\,k}(x)$. An
important object of mathematical study is the zero set of the polynomials $\phi_k(t)$. The following result, which
is no doubt well known, but proved here for completeness, asserts that one can recover the orthogonal polynomial
$\phi_{m}(t),\,m\in I_{n}$, from the characteristic polynomial of $x_m$. In particular $E_x(m) = \{\mu_{1\,m},
\ldots,\mu_{m\,m}\}$
is the zero set of the polynomial $\phi_{m}(t)$. \vs {\bf Theorem 2.20. } {\it Let the notation be as above. In
particular let
$x\in M(n)$ be given by (2.68). For $m\in I_n$ let $\phi_{m}^{(1)}(t) = {1\over a_{m}}\,\phi_{m}(t)$ so
that $\phi_{m}^{(1)}(t) $ is the monic polynomial corresponding to $\phi_{m}(t)$. Then $\phi_{m}^{(1)}(\lambda) $
is the characteristic polynomial of $x_m$ so that $$E_x(m)\,\,\hbox{ is the set of zeros of
${\phi_{m}}(t)$}\eqno (2.70)$$ and one has the interlacing (2.69).} \vs {\bf Proof.} Let $m\in I_n$ so that, 
by (2.67), $x_m$ is the matrix of 
$\Lambda_m\,\tilde {t}| W_m$ with respect to the basis $\{\phi_{j-1}\},\,j\in I_m$, of $W_m$. Let $y \in M(m)$ be
the matrix of $\Lambda_m\,\tilde {t}| W_m$ with respect to the basis $\{t^{j-1}\},\,j\in I_m$, of $W_m$. 
If $j\in I_{m-1}$ then of course $\Lambda_m\,\tilde {t}(t^{j-1}) = t^j\in W_{m}$. But now let $b_i\in \Bbb R$ be
such that 
$$\phi_{m}^{(1)}(t) = t^m + \sum_{j\in I_m} b_j \,t^{j-1}\eqno (2.71)$$ But $\Lambda_m(\phi_{m}^{(1)}(t))= 0$ by
orthogonality. Thus $$\eqalign{ \Lambda_m\,\tilde {t}(t^{m-1})&= \Lambda_m(t^m)\cr &= - \sum_{j\in I_m} b_j
\,t^{j-1}\cr}\eqno (2.72)$$ Hence $y$ is the companion matrix $$\left(\matrix{0&0&\cdots &0&-b_1
\cr  1&0 &\cdots &0&-b_2\cr 
0&1&\cdots&0&-b_3\cr
\vdots&\vdots &\ddots & \vdots&\vdots \cr 0&0&\cdots &1&-b_m\cr}\right)$$ But then $\phi_{m}^{(1)}(\lambda)$ is
the characteristic polynomial of $y$. Consequently $\phi_{m}^{(1)}(\lambda)$ is also the characteristic polynomial
of $x_m$. QED\vs Applying Theorem 2.19 we can say something about the uniqueness of $x\in M(n)$ (see (2.68)) in
yielding the orthogonal polynomials $\phi_m(t),\,m\in I_n$. \vs {\bf Theorem 2.21.} {\it Let $\{\phi_m(t)\}, m\in
I_n$, be the set of normalized orthogonal polynomials, where $deg\,\,\phi_m(t) = m$, for a general measure $\nu$
on $\Bbb R$ specified as above. Then the matrix
$x\in M(n)$, given by (2.68), is the unique symmetric Jacobi matrix with positive adjacent diagonals
 such that
$E_x(m)$, for any $m\in I_n$, as a set, is the set of zeros of the polynomial $\phi_m(t)$.}\vs {\bf Proof.} If $y\in M(n)$
is a Jacobi matrix such that $E_y(m) = E_x(m)$ for any $m\in I_n$ then, by Theorem 2.19, $y$ is a
conjugate of
$x$ by a diagonal matrix $g\in Gl(n)$. But if conjugation by $g$ is to preserve symmetry and the positivity
condition it is immediate that $g$ must be a constant matrix. In such a case of course
$x=y$. QED \vskip 1pc

\rm
 \centerline{\bf 3. The group $A$ and the integration of $\xi_p$
for any $p\in J(n)$ (see (2.30))}\vskip 1pc 3.1. Let $m\in I_n$.
Recall that $P(m)^{Gl(m)}$ (see (2.30)) is the polynomial ring with the $m$ generators $f_{(k,m)},\,k\in I_m$ (see
(2.38)). Let
$\a(m)$ be the commutative (see (2.32)) Lie
algebra of analytic vector fields on $M(n)$ spanned by
$\xi_{(k,m)}$ for $k\in I_m$. Here we are retaining the notation (2.43). 

If $m= n$ then $\a(n) = 0$ by (2.7). Now assume $m\in I_{n-1}$. Note then that $dim\,\a(m) = m$ by (2.43), (2.46),
(2.49) and the existence of strongly regular elements. Let $A(m)$ be a simply
connected complex analytic group where $\a(m) = Lie\,A(m)$. Of course $A(m)\cong \Bbb C^m$. We wish to prove that
$\a(m)$ integrates to an analytic action of $A(m)$ on $M(n)$.

For any
$x\in M(n)$ we recall the subspace $V_x(m)\s T_x(M(n))$ defined as in
(2.36). Since $\{f_{(k,m)}\},\,j\in I_m,$ generate $P(m)^{Gl(m)}$, $$V_x(m),\,\, \hbox{for any
$x\in M(n)$, is spanned by $\{(\xi_{(k,m)})_x\},\,\,k\in I_m$}\eqno (3.1)$$
Consequently if $x\in M(n)$, $$dim\,V_x(m)\leq m\eqno (3.2)$$ However $$dim\,V_x(m) =
m\,\,\hbox{if $x$ is strongly regular, by Theorem 2.7}$$ 

For any $w\in M(m)$ clearly $$ w + M(m)^{\perp} = \{x\in M(n)\mid x_m= w\} \eqno (3.3)$$ using
the notation of (2.37). But obviously one has the fibration $$M(n) = \sqcup_{w\in M(m)} \,\,\,\, w +
M(m)^{\perp}\eqno (3.4)$$\vskip .5pc {\bf Lemma 3.1.} {\it For any $\xi\in \a$ and any $w\in M(m)$
the vector field $\xi$ is tangent to the fiber $w+ M(m)^{\perp}$. That is, 
$$V_x(m)\s T_x(x_m + M(m)^{\perp})\eqno (3.5)$$ for any $x\in M(n)$ (noting of course that $x\in
x_m + M(m)^{\perp}$). } 

\vs {\bf Proof.} Obviously $w + M(m)^{\perp}$ is the nonsingular
affine subvarity of $M(n)$ defined by the equations
$\alpha_{ij}-\alpha_{i\,j}(w) = 0$ for all $\{i,j\}\s I_m$. A vector
field $\eta$ on $M(n)$ is then clearly tangent to $w + M(m)^{\perp}$
if
$\eta(\alpha_{i\,j})= 0$ for all $\{i,j\}\s I_m $. But this is certainly the case if $\eta=\xi$.
Indeed $f_{(k,m)}\in P_(m)^{Gl(m)}$, for $k\in I_m$, and $\alpha_{ij}\in P(m)$. Hence $f_{(k,m)}$
Poisson commutes with $\alpha_{ij}$. QED\vs Recalling the notation of (2.46) let $x\in M(n)$ and let
$G_{x,m}\s Gl(m)$ be the (commutative) subgroup of all $g\in Gl(m)$ such that $g_m$ is an
invertible element in $Z_{x,m}$ (where $Id_m$ is taken as the identity). It is clear that
$G_{x,m}$ is an algebraic subgroup of
$Gl(n)$ and
$Z_{x,m}$ (as a Lie algebra) is the Lie algebra of $G_{x,m}$. Since
$A(m)$ is simply-connected there clearly exists a homomorphism of commutative complex analytic
groups
$$\rho_{x,m}:A(m)\to G_{x,m}\eqno (3.6)$$ whose differential $(\rho_{x,m})_*$ is given by
$$(\rho_{x,m})_*(\xi_{(k,m)}) = -(x_m)^{m-k}\eqno (3.7)$$ for any $k\in I_m$. \vs {\bf Remark 3.2.} Note that 
$G_{x,m}$ is connected since first of all $G_{x,m}\cong \{g_m\mid g\in G_{x,m}\}$ and 
clearly $$\{g_m\mid g\in G_{x,m}\} = \{y\in Z_{x,m}\mid f_{1,m}(y)\neq 0\}$$ It follows then that
$$\rho_{x,m}:A(m)\to G_{x,m},\qquad \hbox{is surjective}\eqno (3.8)$$ since, clearly, $(\rho_{x,m})_*$ is a
surjective Lie algebra homomorphism of 
$\a(m)$ to $Z_{x,m}$.\vs Let
$x\in M(n)$. Now obviously
$M(m)$ and $M(m)^{\perp}$ are both stable under $ad\,\,M(m)$ and $Ad\,\,Gl(m)$. But obviously $x_m$ is fixed by 
$Ad\,G_{x,m}$. Consequently $$x_m + M(m)^{\perp}\,\,\hbox{ is stable under $Ad\,G_{x,m}$}\eqno (3.9)$$\vskip .5pc
{\bf Theorem 3.3.} {\it Let $m\in I_{n-1}$. Then $\a(m)$ integrates to a (complex analytic) action of $A(m)$ on
$M(n)$,
$$A(m)\times M(n)\to M(n),\qquad (a,y)\mapsto a\cdot y\eqno (3.10)$$ More explicitly if $x\in M(n)$ and
$y\in x_m + M(m)^{\perp}$ (noting that any $y\in M(n)$ is of this form by (3.3) and (3.4)) one has $$a\cdot y =
Ad\,(\rho_{x,m}(a))(y)\eqno (3.11)$$ for any $a\in A(m)$.} \vs {\bf Proof.} For the first statement it suffices,
by Lemma 3.1, to show that $\a(m)|(x_m + M(m)^{\perp})$ integrates to an action of $A(m)$ on $x_m +
M(m)^{\perp}$.

Using the notation of \S 2.1 let $Y\s P^1(n)$ be the span of $\alpha_{ij}$ for $(i,j)\in L_m$ so that $dim\,Y =
n^2-m^2$ and $Y$ is the orthocomplement of $M(m)$ in $P^1(n)$. Let ${\cal Y}$ be the affine algebra of functions on
$x_m + M(m)^{\perp}$. For any
$\alpha\in Y$ let
$\alpha_o\in {\cal Y}$ be defined so that if $u\in M(m)^{\perp}$ then $$\alpha_o(x_m + u) = \alpha(u)\eqno (3.12)$$
Let $Y_o =\{\alpha_o\mid \alpha\in Y\}$ so that clearly the subspace $Y_o$ of ${\cal Y}$ generates ${\cal Y}$. Let
$$r:Y\to Y_o\eqno (3.13)$$ be the linear isomorphism defined by putting $r(\alpha) = \alpha_o$. Since any
$\alpha\in Y$ vanishes on $x_m$ note that $$r(\alpha) = \alpha|(x_m + M(m)^{\perp})\,\,\,\,\,\forall \alpha\in Y\eqno
(3.14)$$ Since
$Z_{x,m}\s M(m)$ it is clear that $Y$ is stable under the coadjoint representation of $Z_{x,m}$ and $G_{x,m}$ (resp.
$coad$ and $Coad$) on $P^1(n)$. We now wish to prove that $Y_o$ is stable under $\a(m)$ and in fact, on $Y$, one has
$$\xi\circ r = r \circ coad\,((\rho_{x,m})_*(\xi))\eqno (3.15)$$ for any $\xi\in \a(m)$. To prove (3.15) we first
note that if
$k\in I_m$ and $y\in x_m + M(m)^{\perp}$ then $$(\xi_{(k,m)})_y = (\partial^{[(x_m)^{m-k},y]})_y\eqno (3.16)$$ Indeed 
(3.16) follows from Theorem 2.12 where
we have replaced $x$ by $y$ in (2.43) and recognize that $x_m = y_m$ in our present notation. Applying both
sides of (3.16) to $\alpha\in Y$ one has, by (3.14), $$\eqalign{(\xi_{(k,m)}\alpha_o)(y)
&=(\xi_{(k,m)}\alpha)(y)\cr &=
\langle \alpha,[(x_m)^{m-k},y]\rangle \cr  &= (-coad\,\,(x_m)^{m-k}(\alpha))(y)\cr}$$  But this proves (3.15) (and
hence also the stability of $Y_o$ under $\a(m)$) since $(\rho_{x,m})_*(\xi_{(k,m)}) = -(x_m)^{m-k}$ for $k\in I_m$. 

But now the action of $\a(m)$ on $Y_o$ is a linear representation of $\a(m)$ on a finite dimensional vector space.
Hence it exponentiates to a representation $\pi$ of $A(m)$ on $Y$. But then the commutative diagram (3.15) yields
the commutative diagram $$\pi(a)\circ r = r\circ Coad\,\,(\rho_{x,m}(a))\eqno (3.17)$$ on $Y$ for any $a\in A(m)$.
But now if one defines an action of $A(m)$ on $x_m + M(m)^{\perp}$ by putting $a\cdot y = Ad\,(\rho_{x,m}(a))(y)$ we
observe that, for any $\alpha\in Y$, $$\alpha_o(a^{-1}\cdot y) = (\pi(a)(\alpha_o))(y)\eqno (3.18)$$ Indeed this
follows from (3.17) since $$ \alpha (Ad\,\,\rho_{x,m}(a^{-1})(y)) = (Coad\,\,(\rho_{x,m}(a))(\alpha))(y)$$ But now
replacing
$a$ in (3.18) by $exp\,\,\,t\,\xi$ for $\xi\in \a(m)$ and differentiating it follows from the definition of $\pi$
that the action of $A(m)$ on $x_m + M(m)^{\perp}$ given by (3.11) integrates the Lie algebra $\a(m)|(x_m +
M(m)^{\perp})$. QED\vs 3.2. By Theorem 2.7, and the existence of strongly regular elements (see e.g. (2.27)) the sum
$\a$ of the
$\a(m)$ for $m\in I_{n}$ is a direct sum of the $\a(m)$ where $m\in I_{n-1}$ so that we can write $$\a
=\oplus_{m\in I_{n-1}}\,\,\a(m)\eqno (3.19)$$  We note then that $\a$ is a commutative Lie algebra of analytic
vector fields on $M(n)$ of dimension
$d(n-1)$ with basis $\xi_{p_{(i)}}\,\,i\in d(n-1)$, where for $m\in I_n$, and $k\in I_m$, $$p_{(d(m) +k)} =
f_{(k,m)}\eqno (3.20)$$ (see (2.38)). Let $A = A(1)\times \cdots \times A(n-1)$ so that we can regard $\a= Lie\,
A$ and note that as an analytic group $$A \cong \Bbb C^{d(n-1)}\eqno (3.21)$$ As a consequence of Theorem 3.3
one has
\vs {\bf Theorem 3.4.} {\it The Lie algebra $\a$ of vector fields on $M(n)$ integrates to an action of $A$ on $M(n)$
where if $a\in A$ and $a = (a(1),\ldots, a(n-1))$ for $a(m)\in A(m)$ then $$a\cdot y =
a(1)\cdot(\cdots (a(n-1)\cdot y)\cdots )\eqno (3.22)$$ for any $y\in M(n)$. Furthermore the orbits of $A$ on
$M(n)$ are leaves of the distribution, $x
\mapsto V_x$, on $M(n)$. See 
\S 2.4, Remark 2.9 and Theorem 2.13.}\vs {\bf Proof.} Since $\a$ is commutative the actions of $A(m)$ for all
$m\in I_{n-1}$ clearly commute with one another. Thus (3.22) defines an action of $A$ on $M(n)$ which, by Theorem 3.3,
integrates
$\a$. The final statement of the theorem follows from Theorem 2.12 and (2.48). QED\vs Recall that $J(n)\s P(n)$ is
the Poissson commutative associative subalgebra defined in \S 2.4. See (2.30). 
\vs {\bf Theorem 3.5.} {\it Let $p\in J(n)$ be arbitrary. Then the vector field $\xi_p$ is globally integrable defining
an action of $\Bbb C$ on $M(n)$. Furthermore this action commutes with $A$ and stabilizes any $A$--orbit. 

Let $\{p_i'\},\,i\in I_d$, for some integer $d$, be an arbitrary set of generators of $J(n)$. Let $\a'$ be the
(commutative) Lie algebra spanned by $\xi_{p_i'},\,i\in I_d$, and let $A'$ be a corresponding simply connected Lie
group. Then $\a'$ integrates to an action of $A'$ on $M(n)$. Furthermore the action of $A'$ commutes with action of
$A$ and the orbits of
$A'$ are the same as the orbits of $A$. In particular any $A'$-orbit is a leaf of the distribution $x\mapsto V_x$ on
$M(n)$.}\vs {\bf Proof.} Let
$y\in M(n)$ and let
$Q$ be the orbit
$A\cdot y$ of
$A$ so that
$Q$ is an analytic submanifold of $M(n)$. By Theorem 3.4 one has $$V_x = T_x(Q)\,\,\,\,\forall x\in Q\eqno
(3.23)$$ Let
$\p = \{f\in \a\mid \,\,\xi_f|Q = 0\}$ and let $\q$ be a linear complement of $\p$ in $\q$. Then by the
commutativity of $A$ it follows that if $dim\,Q = \ell$ then $\dim\,\q = \ell$ and if $q_j,\,j\in I_{\ell}$, is a
basis of $\q$ then $\{(\xi_{q_j})_x\},\,j\in I_{\ell}$, is a basis of $T_x(Q)$ for all $x\in Q$. But now 
$\xi_p| Q$ is tangent to $Q$ by (3.23). Thus there exists functions $h_j,\,j\in I_{\ell}$ on $Q$ so that $$\xi_p =
 \sum_{j\in I_{\ell}}\,\,h_j \xi_{q_j}\eqno (3.24)$$ on $Q$. But the commutativity $[\xi_p, \xi_{q_j}] = 0$ implies
that the functions $h_j$ are constant on $Q$. By Theorem 3.3 the vector field $\xi_{q_j}|Q$ integrates to an action,
$t\mapsto exp\,\,t\,\xi_{q_j}|Q$ of $\Bbb C$ on $Q$. But then one gets another action of $\Bbb C$ on $Q$ by putting 
$$b(t) = exp\,\,t\,h_1\,\xi_{q_1}\cdots exp\,\,t\,h_{\ell}\,\xi_{q_{\ell}}|Q$$ But then it is immediate from (3.24)
that $\xi_p|Q$ integrates to the action $b(t)$ on $Q$. Thus $\xi_p$ integrates to an action of $\Bbb C$ on
$M(n)$. Clearly this action of
commutes with the action of $A$ on $M(n)$ by (2.32). In addition the argument above implies that any $A$--orbit
is stable under the action of $\Bbb C$. 

It then follows that $\a'$ integrates to an action of $A'$ on $M(n)$ and the action of $A'$ commutes with the
action of $A$. In addition the argument above implies that that any
$A$--orbit
$Q$ is stable under the action of $A'$. But since the $p_i'$ generate $J(n)$ the span of $(dp_i')_x$, for $i\in I_d$,
is the same as the span of $(dp_{(j)})_x$, for $j\in I_{d(n)}$, for any $x\in M(n)$. Thus $(\xi_{p_i'})_x$, for $i\in
I_d$, spans $V_x$ by (1.14). Consequently, all orbits of $A^{\prime}$ on $Q$ are
open. But $Q$ is connected since $A$ is connected. 
  But then $A^\prime$ is obviously transitive on $Q$.  QED

\vs 3.3. In this section we will give some properties of the group $A$. We first
observe the elementary property that
$A$ operates ``vertically" with respect to the ``fibration" (2.9) of $M(n)$. Recall (see (2.8)) the morphism
$\Phi_n:M(n)
\to \Bbb C^{d(n)}$. Regarding 
$c_j$ as a coordinate function on $\Bbb C^{d(n)}$ note that, by (2.8), for $j\in I_{d(n)}$,  $$c_j\circ \Phi_n =
p_j
\eqno (3.25)$$ \vskip .5pc {\bf Proposition 3.6.} {\it The group $A$ stabilizes the ``fiber" $M_c(n)$ of $M(n)$
(see (2.9)) for any $c\in \Bbb C^{d(n)}$. In fact $M_c(n)$ is stabilized by the action of $\Bbb C$ defined in Theorem
3.5.}
\vs {\bf Proof.} For
$i,j\in I_{d(n)}$ one has
$\xi_{p_{(i)}}p_j = 0$ by (2.31) (see (3.20)). But then the function $p_j$ on $M(n)$ is invariant under $A$. Hence
$M_c(n)$ is stabilized by
$A$, for any $c\in \Bbb C^{d(n)}$, by (3.25). The final statement then follows from Theorem 3.5. QED\vs 
Let
$x\in M(n)$. We wish to give an explicit description of the $A$-orbit $A\cdot x$. Obviously $G_{x,m}$ is an algebraic
subgroup of
$Gl(n)$, for any 
$m\in I_n$,
and the adjoint representation of $G_{x,m}$ on $M(n)$ is an algebraic representation. In the following theorem we
use the term ``constructible set". This is a concept in algebraic geometry due to C. Chevalley. For its precise
definition and properties see
\S 5, Chapter II in [C] or Exercise 1.9.6, p. 19 in [S], or \S 1.3, p. 3 in [B]. \vs {\bf Theorem 3.7.} {\it Let
$x\in M(n)$. Consider the following morphism of nonsingular irreducible (see Remark 3.2) affine varieties
$$G_{x,1}\times
\cdots \times G_{x,n-1}
\to M(n)\eqno (3.26)$$ where for
$g(m)\in G_{x,m},\,\,m\in I_{n-1}$, $$ (g(1),\ldots, g(n-1))\mapsto Ad\,(g(1))\cdots Ad\,(g(n-1))(x)\eqno (3.27)$$
Then the image of (3.26) is exactly the
$A$-orbit $A\cdot x$. In fact if $a\in A$ and $a= a(1)\cdots a(n-1)$ where $a(m)\in A(m)$ then $a\cdot x$ is given by
the right side of (3.27) where $g(m) = \rho_{x,m}(a(m))$. In particular
$A
\cdot x$ is an irreducible, constructible (in the sense of Chevalley) subset of $M(n)$. Furthermore the Zariski closure
of
$A\cdot x$ is the same as its usual Hausdorf closure. In addition $A\cdot x$ contains a Zariski open subset of its
closure.}
\vs {\bf Proof.} Let $m\in I_{n-1}$. For $g(m)\in G_{x,m}$ there exists, by (3.8) and Theorem 3.2, $a(m)\in A(m)$
such that
$a(m) = Ad\,g(m)$ on $x_m + M(m)^{\perp}$. Let
$b(m) = a(m)\cdots a(n-1)$ and $h(m) = g(m)\cdots g(n-1)$. Obviously the right side of (3.27) can be written
$Ad\,h(1)(x)$. Inductively (downwards) assume that $b(m+1)\cdot x = Ad\,h(m+1)(x)$ for
$m<n-1$. Note that the induction assumption is satisfied if $m+1 = n-1$ by Theorem 3.2. But since $p_{(i)}$ is
invariant under
$Gl(m)$ for $i>d(m-1)$ note (key observation) that
$(b(m+1)\cdot x)_m = x_m$ so that 
$b(m+1)\cdot x\in x_m + M(m)^{\perp}$. Thus $b(m)\cdot x = Ad\,h(m)(x)$ by Theorem 3.2. Hence, by
induction, $Ad\,h(1)(x)\in  A\cdot x$. Conversely let $a\in A$ so that we can write $a = a(1)\cdots a(n-1)$ for 
$a(m)\in A(m)$. Let $g(m) =
\rho_{x,m}(a(m))$ so that $g(m)\in G_{x,m}$. Then the argument above establishes that $a\cdot x= Ad\,h(1)(x)$.
 Hence the image of (3.26) is $A\cdot x$. The image of an irreducible set under a morphism is
obviously irreducible. But the image $A\cdot x$ is constructible by Proposition 5 in Chapter V, p. 95-96 in [C] or
by Proposition, p. 4 in [B]. See also Corollary 2, \S 8, p.51 in [M]. But then the Euclidean closure of $A\cdot x$ is
the same as its Zariski closure by Corollary 1, \S 9, p.60 in [M]. On the other hand the irreducible
constructible set
$A\cdot x$ contains a Zariski open subset of its closure by Proposition 4 in Chapter 5, p.95 in [C] or by
Proposition, p. 4 in [B]. QED\vs {\bf Remark 3.8.} In (3.22) the ordering of the terms $a(m)$ is immaterial since
these elements commute with one another. However the ordering of the terms $g(m)$ in (3.27) cannot be permuted in
general. The groups $G_{x,m}$ do not, in general, commute with one another. The point is that $A(m)$ is given by
$Ad\,G_{x,m}$ only on $x_m + M(m)^{\perp}$. \vs By the final statement of Theorem 3.4 and (2.34) one has, for any
$x\in M(n)$, 
$$dim\,A\cdot x \leq d(n-1)\eqno (3.28)$$ and, by (2.35), $$dim\,A\cdot x = d(n-1) \iff \hbox{$x$ is strongly
regular (see \S 2.3)}\eqno (3.29)$$ That is, (see \S 2.3) the nonempty Zariski open subset $M^{sreg}(n)$ of $M(n)$
is stable under the action of
$A$ and $$\hbox{the $A$-orbits in $M^{sreg}(n)$ are exactly all the $A$-orbits of maximal dimension ($d(n-1)$)}\eqno
(3.30)$$ We will now see that we can make a much stronger statement about $A\cdot x$ than that in Theorem 3.6 in
case $x$ is strongly regular. We first need some preliminary results. 

Let $c\in \Bbb C^{d(n)}$. One has $M_c(n) = (\Phi_n)^{-1}(c)$ (see (2.8)) so that $M_c(n)$ is a Zariski closed
subset of $M(n)$. Now let 
$M_c^{sreg}(n) = M_c(n)\cap M^{sreg}(n)$. \vs {\bf Remark 3.9.}. Note that $M_c^{sreg}(n)$ is a Zariski open
subset of $M_c(n)$. Furthermore $M_c^{sreg}(n)$ is not empty since
$M_c(n)\cap (e+\b)$ is a one point subset of $M_c^{sreg}(n)$, by Theorem 2.5 and (2.27). However $M_c^{sreg}(n)$ may
not be dense in $M_c(n)$. In fact one can show that it is not dense if $n=3$ and $c=0$.\vs  Let $N(c)$ be the number
of irreducible components of
$M_c^{sreg}(n)$ and let
$M_{c,i}^{sreg}(n),\,i=1,\ldots,N(c)$, be some indexing of these components. In the following proposition overline
 is
Zariski closure. 
\vs {\bf Proposition 3.10.} {\it The irreducible component decomposition of $\overline {M_c^{sreg}(n)}$ is
$$\overline {M_c^{sreg}(n)} = \bigcup_{i=1}^{N(c)}\overline  {M_{c,i}^{sreg}(n)}\eqno (3.31)$$ noting that
$\overline  {M_{c,i}^{sreg}(n)}$ is a closed subvariety of $M_c(n)$. Furthemore one recovers $M_{c,i}^{sreg}(n)$
from its closure by
$$M_{c,i}^{sreg}(n) = (\overline  {M_{c,i}^{sreg}(n)})\cap M^{sreg}(n)\eqno (3.32)$$ so that $M_{c,i}^{sreg}(n)$
is an Zariski open subvariety of its Zariski closure. Furthermore $M_{c,i}^{sreg}(n)$ is a constructible set so that
its Zariski closure is the same as its usual Hausdorf closure (see p. 60 in [M])} \vs  {\bf Proof.}  The 
equality (3.31)
 is immediate. The statement that (3.31) is an irreducible component decomposition is an easy consequence of
(3.32). In fact it is a general result. See the last paragraph in \S 1.1, p.3 in [B]. Clearly the left side of 
(3.32)
is contained in the right hand side. But the right hand side is an irreducible set in $M^{sreg}_c(n)$. Thus one has
(3.32) by the irreducible maximality of $M_{c,i}^{sreg}(n)$ in $M_c^{sreg}(n)$. We have used Proposition 2 in Chapter
II of [C], p.35 throughout. $M_{c,i}^{sreg}(n)$ is a constructible set by (3.32) and
Proposition 3 in Chapter II of [C], p. 94. The last statement of the proposition follows
from Corollary 1, \S 9, p. 60 in in [M]. QED\vs {\bf Proposition 3.11}. {\it
Let the assumptions and notations be as in Proposition 3.10. Then $M_{c,i}^{sreg}(n)$ and
hence  $\overline  {M_{c,i}^{sreg}(n)}$ have dimension $d(n-1)$ for all $i$. Furthermore
$M_{c,i}^{sreg}(n)$ is nonsingular in $\overline{M_c^{sreg}(n)}$ and a fortiori
$M_{c,i}^{sreg}(n)$ is nonsingular in $\overline  {M_{c,i}^{sreg}(n)}$.}\vs {\bf Proof}.
The proof is an application of Theorem 4 in \S 4 of Chapter 3 in [M], p. 172, where in the
notation of that reference, $X= M(n)$, $Y = \overline{M_c^{sreg}(n)}$, $U$ runs through
all open affine subvarieties of a finite affine cover of $M^{sreg}(n)$, $k =
d(n)$, $f_i = p_i-c_i$. QED \vs Let $c\in \Bbb C^{d(n)}$ be arbitrary. Then since $M_c(n)$ is stable under the action
of $A$ it follows from (3.30) that $M^{sreg}_c(n)$ is stable under $A$.  \vs {\bf Theorem 3.12.} {\it Let the
assumptions and notations be as in Proposition 3.10. Then there are exactly $N(c)$ orbits of $A$ in $M_c^{sreg}(n)$
and these orbits are identical to the irreducible components $M_{c,i}^{sreg}(n)$ of $M_c^{sreg}(n)$. In particular
all the orbits are nonsingular algebraic varieties of dimension $d(n-1)$. } \vs {\bf Proof.} Let $x\in M^{sreg}_c(n)$
so that
$A\cdot x
\s M^{sreg}_c(n)$. But by irreducibility of $A\cdot x$ (see Theorem 3.6) the Zariski closure $\overline {A\cdot x}$
is irreducible. Thus, by Proposition 3.10, there exists $i$ such that $$\overline {A\cdot x}\s \overline
{M_{c,i}^{sreg}(n)}\eqno (3.33)$$ But then if $k =
dim\,\overline {A\cdot x}$ one has $k\leq d(n-1)$ and one has equality in (3.33) if
$k=d(n-1)$. But $A\cdot x$ contains a Zariski open set $U$ of $\overline {A\cdot x}$ by the last statement in Theorem
3.7. Let
$W\s U$ be the nonsingular Zariski open subvariety of simple points in $U$. Then $W$ is an analytic submanifold of
dimension
$k$ in
$M(n)$. But, as a nonempty Zariski open subvariety of a closed subvariety,  
$W = V\cap \overline {A\cdot x}$ where $V$ is a Zariski open subset of $M(n)$. Thus $W$ is open in $A\cdot
x$. But $A\cdot x$ is a $d(n-1)$-dimensional analytic submanifold of $M(n)$. Thus $k = d(n-1)$, by invariance of
dimension of submanifolds, so that 
$$\overline {A\cdot x} = 
\overline {M_{c,i}^{sreg}(n)}\eqno (3.34)$$ But now (3.34) implies that $\overline {M_{c,i}^{sreg}(n)}$ is stable
under the action of $A$ since $\overline {A\cdot x}$ is the closure of $A\cdot x$ in the usual Euclidean topology. 
But then $M_{c,i}^{sreg}(n)$ must be stable under the
action of
$A$ by (3.32). In particular $A\cdot x \s M_{c,i}^{sreg}(n)$ since, clearly, $x\in M_{c,i}^{sreg}(n)$ by (3.32) and
(3.33). Assume
$A\cdot x
\neq M_{c,i}^{sreg}(n)$ and let
$y\in  M_{c,i}^{sreg}(n)-A\cdot x$. Then clearly we may replace $x$ by $y$ in (3.34) so that $\overline {A\cdot y} = 
\overline {M_{c,i}^{sreg}(n)}$. But then $A\cdot x $ and $A\cdot y$ contain disjoint nonempty Zariski open subsets of
$\overline {M_{c,i}^{sreg}(n)}$. This contradicts the irreducibility of $\overline {M_{c,i}^{sreg}(n)}$. Thus we have
proved that $M_{c,i}^{sreg}(n)$ is an $A$-orbit. If $j\neq i$ we may find $z\in M_{c,j}^{sreg}(n)$ which does not 
lie
in any other component of $M_c^{sreg}(n)$. The argument above then implies that $M_{c,j}^{sreg}(n) = A\cdot z$.
QED\vs 3.4.
 Let
$x\in M^{sreg}(n)$. Then by Theorem 3.6 and Theorem 3.12 one has the following surjective morphism of nonsingular
varieties.
$$\nu_x:G_{x,1}\times \cdots \times G_{x,n-1}\to A\cdot x\eqno (3.35)$$ where $$\nu_x(g(1),\ldots,g(n-1)) =
Ad(g(1))\cdots Ad(g(n-1))(x)\eqno (3.36)$$ We now prove \vs {\bf Lemma 3.13.} {\it The map (3.35) is injective so
that in fact $\nu_x$ is bijective for any $x\in M^{sreg}(n)$.} \vs {\bf Proof.} Assume $g(m),g(m)'\in G_{x,m}$ for
$m=1,\ldots,n-1$, and
$$Ad(g(1))\cdots Ad(g(n-1))(x)= Ad(g(1)')\cdots Ad(g(n-1)')(x)$$ By upward induction we will prove that
$g(m)=g(m)'$ for all $m$. The inductive argument to follow also establishes that the result is true for $m=1$. Assume
we have proved that $g(i) = g(i)'$ for $i<m$. Then if $y= Ad(g(m))\cdots Ad(g(n-1))(x)$ one has $y = Ad(g(m)')\cdots
Ad(g(n-1)')(x)$. Let $z=  Ad(g(m))^{-1}(y)$ and $z'=Ad(g(m)')^{-1}(y)$. Then $z_{m+1} = z'_{m+1} = x_{m+1}$. But
$z = Ad(h(m))(z')$ where $h(m) = g(m)^{-1}g(m)'$. But then $h(m)$ commutes with $x_{m+1}$. Thus $h(m)-Id_n\in
Z_{x,m}\cap Z_{x,m+1}$ by (2.46). But $Z_{x,m}\cap Z_{x,m+1} = 0$ by (2.49). Thus $g(m) = g(m)'$ and hence $\nu_x$ is
injective. QED\vs 
We can now describe the structure of all $A$-orbits of maximal (i.e., $d(n-1)$ ) dimension. \vs
{\bf Theorem 3.14.} {\it Let $x\in M^{sreg}(n)$ and for $m=1,\ldots,n-1$ let $G_{x,m}$ be the centralizer of $x_m$ in
$Gl(m)$ so that $G_{x,m}$ is a connected commutative algebraic group of dimension $m$. Then the morphism (3.35) is
an algebraic isomorphism of nonsingular varieties so that as a variety $$A\cdot x \cong G_{x,1}\times \cdots \times
G_{x,n-1}\eqno (3.37)$$}\vs {\bf Proof.} Since $\nu_x$ is a bijective morphism of varieties of dimension $d(n-1)$ it
follows that $\nu_x$ is birational by Theorem 3 of \S IV in [C], p. 115 ($n=0$ in the notation of this reference). 
But then $\nu_x$ is an isomorphism of algebraic varieties by Theorem, p.78, \S 18 of Chapter AG in [B] and also
Theorem 5.2.8, p. 85 in [S] since the range and domain of
$\nu_x$ are nonsingular by Theorem 3.12. QED \vs {\bf Example 3.15}. Since $e\in e + \b$ note that the principal
nilpotent element
$e$ is strongly regular. If
$x= e$ then one readily sees that $A\cdot x$ is the set of all principal nilpotent elements in the nilpotent Lie algebra
$\u'$ of all strictly lower triangular matrices and $\overline {A\cdot x} = \u'$. Writing the set of all principal
nilpotent elements in $\u'$ uniquely, using the right--hand side of (3.36), when $x = e$,  may be new.
vs 3.5. In this section we will assume $x$ satisfies the eigenvalue disjointness condition. That is $x\in
M_{\Omega}(n)$. See \S 2.4 and more specifically (2.53). We recall that $M_{\Omega}(n)$ is a Zariski open subset of
$M^{sreg}(n)$ so that, in particular, Theorem 3.14 applies to $x$. \vs {\bf Theorem 3.16.} {\it Let $x\in M(n)$
satisfy the eigenvalue disjointness condition. Then $Z_{x,m}$ is a Cartan subalgebra of $M(m)$, for any $m\in I_n$,
and $G_{x,m}$ is a maximal (complex) torus in $Gl(m)$ so that $$G_{x,m}\cong (\Bbb C^{\times})^m\eqno (3.38)$$
In addition the orbit $A\cdot x$ is an algebraic subvariety of $M(n)$ and as such $$A\cdot x \cong (\Bbb
C^{\times})^{d(n-1)}\eqno (3.39)$$}
\vs {\bf Proof.}
$Z_{x,m}$ is a Cartan subalgebra of $M(m)$ since $x_m$ is a regular semisimple element of $M(m)$. But this of course
implies that
$G_{x,m}$ is a maximal (complex) torus in $Gl(m)$ since $Lie\,G_{x,m} = Z_{x,m}$. One then has (3.38). Theorem 3.14
then implies the remaining statements of Theorem 3.16.  QED\vs For $m\in I_{n-1}$ let $Y(m)$ be the span of
$e_{i\,m+1},\,i\in I_{m}$ so that $Y(m)$ is an $m$-dimensional subspace of $M(m)^{\perp}$. It is clear that $Y(m)$ is
stable under $Ad\,Gl(m)$ and $ad\,M(m)$ and hence $Y(m)$ inherits, respectively, the structure of a $Gl(m)$ and a
$M(m)$-module. On the other hand $\Bbb C^{m}$ is a module for $Gl(m)$ and $M(m)$ with respect to the natural action
of $Gl(m)$ and $M(m)$. One immediately observes that the linear isomorphism $$Y(m)\to \Bbb C^{m},\,\,\,y\mapsto
(\alpha_{1\,m+1}(y),\ldots,\alpha_{m\,m+1}(y))\eqno (3.40)$$ is an isomorphism of $Gl(m)$ and $M(m)$-modules. 

Let $H(m)$ be a maximal torus of $Gl(m)$ and let $\hh(m) = Lie\,H(m)$ be the corresponding Cartan subalgebra of
$M(m)$. Then, by restriction, $Y(m)$ is a $H(m)$-module with respect to $Ad$ and an $\hh(m)$-module with respect to
$ad$. One notes in fact that $Y(m)$ is a cyclic
$H(m)$-module and a cyclic $\hh(m)$-module. This is
transparent, from the module isomorphism (3.40), in case we choose
$H(m) = Diag(m)$ where $Diag(m)$ is the group of diagonal matrices in $Gl(m)$. The general case follows by
conjugation. An element
$y\in Y(m)$ will be called a cyclic $H(m)$-generator in case it generates $Y(m)$ under the action of $H(m)$. A similar
terminology will be used for $\hh(m)$. By going to the diagonal case it is
immediate that $y$ is a cyclic generator for $H(m)$ if and only if it is a cyclic generator for $\hh(m)$. Indeed
let $$Y^{\times}(m) =
 \{y\in Y(m)\mid \alpha_{i\,m+1}(y)\neq 0,\,\,\forall i\in I_m\}\eqno (3.41)$$ Then clearly $y$ is a cyclic
generator for $Diag(m)$ and a cyclic generator for $Lie\,Diag(m)$ if and only if $y\in Y^{\times}(m)$. 
The following proposition is established by again going to the diagonal case and conjugating (using of course the
module isomorphism (3.40)).\vs  {\bf Proposition 3.17.} {\it Let
$m\in I_{n-1}$. Let
$H(m)$ be a maximal torus of $GL(m)$ and let $\hh(m) = Lie\,H(m)$. Let $y\in Y(m)$. Then the following conditions
are all equivalent. $$\eqalign{&(a)\,\,\,y\,\,\hbox{is a cyclic generator for $H(m)$ (or $\hh(m)$)}\cr
&(b)\,\,\,Ad\,H(m)(y)\,\,\hbox{is the unique $Ad\,H(m)$-orbit in $Y(m)$ of maximal dimension ($m$)}\cr
&(c)\,\,\,Ad\,h(y) = y\,\,\hbox{for $h\in H(m)$ implies $h = Id_n$}\cr
&(d)\,\,\,\hbox{ If $w\in Y(m)$ then }\,\,w \in Ad\,H(m)(y)\,\,\iff w\,\,\hbox{is a cyclic generator for
$H(m)$}\cr &(e)\,\,\,\{[z_{(i)},y]\},\,i\in I_m,\,\,\hbox{is a basis of $Y(m)$ if $\{z_{(i)}\},\,i\in I_m$, is a
basis of
$\hh(m)$}\cr}$$}\vskip .5pc For any $x\in M_{\Omega}$ and $m\in I_{n-1}$ let
$$Y^x(m)= \{y\in Y(m)\mid y\,\,\hbox{ is a cyclic generator for $G_{x,m}$}\}\eqno (3.42)$$

We recall that $x_m$ is the ``cutoff" of $x$ in $M(m)$ for any $x\in M(n)$ where
$m\in I_n$. Henceforth, if $m\in I_{n-1}$, let $x_{\{m\}}\in Y(m)$ be the ``component" of $x$ in $Y(m)$. That is,
$x_{\{m\}} = \sum_{i=1}^m \alpha_{i\,m+1}(x)e_{i\,m+1}$. 
\vs{\bf Proposition 3.18.} {\it Let $x\in M_{\Omega}(n)$. Then $x_{\{m\}}\in Y^{x}(m)$.}\vs  {\bf Proof.} Let $g\in
Gl(m)$ be such that $g\,x_m\, g^{-1}\in Diag(m)$. Then clearly $g\,G_{x,m}\,g^{-1} = Diag(m)$. It suffices then to
show that 
$g\,x_{\{m\}}\, g^{-1}\in Y^{\times}(m)$. To prove this we use the notation and the argument in the proof of
Theorem 2.17. But then the result follows from (2.60). QED\vs Proposition 3.17 lists a number of criteria for $y\in
Y(m)$ to be an element of $Y^x(m)$. It is convenient to add the following to this list. \vs {\bf Proposition 3.19.}
{\it Let $x\in M_{\Omega}(n)$ and let $m\in I_{n-1}$. Let $y\in Y(m)$. Then $y\in Y^x(m)$ if and only if
$\{(ad\,x_m)^j(y)\},\,j = 0,\,\ldots,m-1,$ is a basis of $Y(m)$.}\vs {\bf Proof.} As noted in Proposition 3.17 one
has $y\in Y^x(m)$ if and only if $y$ is a cyclic $Z_{x,m}$ generator. But $x_m$ is a regular semisimple element of
$M(m)$. But then $ad\,x_m|Y(m)$ is diagonalizable with $m$ distinct eigenvalues, by the module isomorphism (3.40). But
this readily establishes the proposition by standard linear algebra. QED \vs {\bf Remark 3.20.} If $x\in
M_{\Omega}(n)$ note that
$Y^x(m)$ is a Zariski open
$Ad\,G_{x,m}$-orbit in
$Y(m)$ by Proposition 3.17. However, if
$m>1$ the maximal torii $G_{x,m}$ of $Gl(m)$ do not run through all maximal torii of $Gl(m)$ as $x$ runs through
$M_{\Omega}(n)$. For example one easily has $G_{x,m}\neq Diag(m)$ for any $x\in M_{\Omega}(n)$. This restricts the
possible sets $Y^x(m)$. In fact, for example, the following theorem implies that $Y^x(m)\neq Y^{\times}(m)$ for all
$x\in M_{\Omega}(n)$. \vs  {\bf Theorem 3.21.} {\it Let $m\in I_{n-1}$. One has $e_{m\,m+1}\in Y^x(m)$ for all
$x\in M_{\Omega}(n)$.} \vs {\bf Proof.} Let $x\in M_{\Omega}(n)$. If $m=1$ the result is obvious. Assume $m>1$. Let $W$
be the subspace of
$Y(m)$ spanned by $(ad\,x_m)^j(e_{m,m+1})$ for $j=0,\ldots,m-1$. Using the characteristic polynomial of $ad\,x_m$ it
is clear that $W$ is stable under $ad\,x_m$. To show that
$e_{m,m+1}$ is a cyclic
$G_{x,m}$ generator ($=$ cyclic $Z_{x,m}$ generator) it suffices, by Proposition 3.19, to show that 
$W= Y(m)$. For this, of course, it suffices to show $$dim\,W = m\eqno (3.43)$$ Let $Y_o(m)$ be the span of
$e_{i\,m+1},\, i=1,\ldots,m-1$, so that
$dim\,Y_o(m) = m-1$. Let $Q:Y(m)\to Y_o(m)$ be the projection with respect to the
decomposition $Y(m) = Y_o(m) \oplus \Bbb C\, e_{m\,m+1}$. Let
$W_o= W\cap Y_o(m)$. Since $e_{m\,m+1}\in W$ one has $W = W_o\oplus \Bbb C\,e_{m\,m+1}$. To prove (3.43) it suffices to
show that $W_o = Y_o(m)$. Clearly $Q(W) = W_o$. Let
$v\in W_o$. Then, since $\alpha_{m\,m+1}(v)= 0$, note that $$Q[x_m,v] = [x_{m-1},v]\eqno (3.44)$$ so that $W_o$
is stable under $ad\,x_{m-1}$. But now 
$[x_m,e_{m\,m+1}]\in W$ and hence $Q[x_m,e_{m\,m+1}]\in W_o$. Thus if $w_j =
(ad\,x_{m-1})^j(Q[x_m,e_{m\,m+1}])$ then 
$$w_j\in W_o\eqno (3.45)$$ for $j\in \{0,\ldots,m-2\}$.
 But note that
$$[Q[x_m,e_{m\,m+1}],e_{m+1\,m}] = x_{\{m-1\}}\eqno (3.46)$$ But $ad\,x_{m-1}$ commutes with
$ad\,e_{m+1\,m}$. Thus $$[w_j,e_{m+1\,m}] = (ad\,x_{m-1})^jx_{\{m-1\}}\eqno (3.47)$$
But $x_{\{m-1\}}\in Y^{x}(m-1)$ by Proposition 3.18. Hence the dimension of the
subspace spanned by the vectors on the right--hand side of (3.47) for
 $j\in \{0,\ldots,m-2\}$ is $m-1$ by Proposition 3.19. Thus the space spanned by the
$w_j$ must have dimension $m-1$. This proves that $W_o = Y_o(m)$. QED \vs For any $x\in M(n)$ let $x^T$ be the
transpose matrix and for any subset $X\s M(n)$ let $X^T=\{x^T\mid x\in X\}$. Obviously $$p_i(x) = p_i(x^T)
\eqno(3.48)$$
for all $x\in M(n)$ and
$i\in I_{d(n)}$ so that $$M_c(n) = M_c(n)^T\eqno (3.49)$$ for any $c\in \Bbb C^{d(n)}$. Note that (see \S 2.1 and
(2.12))
$$-e^T = \sum_{m\in I_{n-1}} e_{m\,m+1}\eqno (3.50)$$

Let $\u\s M(n)$ be the Lie algebra  of strictly upper triangular matrices. Then clearly $$\u = \oplus_{m\in I_{n-1}}
Y(m)\eqno (3.51)$$\vskip .5pc {\bf Lemma 3.22.} {\it Let $x\in M_{\Omega}(n)$. Then there exists $b_x\in A$ such that
$$b_x\cdot x \in (-e + \b)^T \eqno (3.52)$$}\vs {\bf Proof.} Let $m\in I_{n-1}$. Assume inductively (downward) we have
found
$g(k)\in G_{x,k}$ for $k= m,\ldots,n-1$ such that if $h(m) = g(m)\cdots g(n-1)$ and $y(m) = (Ad\,h(m))(x)$ then
$$y(m)_{\{j\}}= e_{j\,j+1}\,\,\,\hbox{for}\,\,j= m,\ldots,n-1$$ If $m=n-1$ this induction assumption is satisfied, by 
Theorem 3.21, in that we can choose $g(n-1)\in G_{x,n-1}$ such that $Ad\,g(n-1)(x_{\{n-1\}}) = e_{n-1\,n}$. To
advance the induction assume $m>1$ and note that, since $x_k$ is fixed by $Ad\,G_{x,k}$ for all $k$, one has $(y(m))_m
= x_m$. Hence
$y(m)_{\{m-1\}} = x_{\{m-1\}}$. By Theorem 3.21 we can choose $g(m-1)\in G_{x,m-1}$ so that
 $$Ad\,g(m-1)(x_{\{m-1\}}) = e_{m-1\,m}\eqno (3.53)$$ Let $h(m-1) = g(m-1)h(m)$ and $y(m-1) =  (Ad\,h(m-1))(x)$. The
point is that the induction assumption is satisfied for $y(m-1)$ since the element $e_{j\,j+1}\in Y(j)$ is clearly
fixed under the adjoint action of $Gl(m-1)$ on $Y(j)$ for $j= m,\ldots,n-1$. By Theorem 3.6 there exists $b_x\in A$ such
that $b_x\cdot x = (Ad\,h(1))(x)$. But then $(b_x\cdot x)_{\{k\}} = e_{k\,k+1}$ for $k= 1,\ldots,n-1$. But then by
(3.50) one has $b_x\cdot x + e^T\in \b^T$. QED\vs Let $c\in \Omega(n)$ where we recall $\Omega(n)\s \Bbb C^{d(n)}$ is
the Zariski open set defined by the eigenvalue disjointness condition. See \S 2.4 and (2.53). We now have the following
neat discription of the fiber $M_c(n)$ of (2.8).\vs {\bf Theorem 3.23.} {\it Let $c\in \Omega(n)$. See (2.53). Then the
group
$A$ operates transitively on the fiber $M_c(n)$ of (2.8). Furthermore $M_c(n)$ is an $A$-orbit of maximal dimension in
$M(n)$. Moreover $M_c(n)$ is a nonsingular Zariski closed subvariety of $M(n)$ of dimension $d(n-1)$. As an
algebraic variety $$M_c(n)\cong (\Bbb C^{\times})^{d(n-1)}\eqno (3.54)$$}\vs {\bf Proof.} Since $M_c(n)\s M^{sreg}(n)$
(see (2.55)) one has $M_c^{sreg}(n) = M_c(n)$ in the notation of Remark 3.9. But obviously $M_c(n)$ is Zariski closed
in $M(n)$. But, by Theorem 3.12, the irreducible components of $M_c(n)$ are then the $A$-orbits in $M_c(n)$, each of
which is a maximal orbit and a nonsingular variety of dimension $d(n-1)$. We now show that there is only one orbit
(i.e., $N(c) = 1$ in the notation of Theorem 3.12). By Remark 2.4 we may replace $e$ by $-e$ in Theorem 2.5 and,
recalling (3.48),
$-e +\b$ by $(-e + \b)^T$ so that the restriction $$\Phi_n:(-e + \b)^T\to \Bbb C^{d(n)}\eqno (3.55)$$ is an algebraic
isomorphism. But now if $x,y \in M_c(n)$ there exists, by Lemma 3.22, $b_x,b_y\in A$ such that both $b_x\cdot x$ and
$b_y\cdot y$ are in $(-e + \b)^T$. But since $M_c(n)$ is stabilized by $A$ it follows from (3.55) that $b_x\cdot x =
b_y\cdot y$. Thus $x$ and $y$ are $A$ conjugate and hence $N(c) = 1$. The isomorphism (3.54) then follows from the
isomorphism (3.39). QED\vs The following two results could have been proved at an earlier point. However they are
manifestly transparent as a consequence of Theorem 3.23.\vs {\bf Proposition 3.24.} {\it Let $c\in \Omega(n)$ and let
$x\in M_c(n)$. Then the tangent vectors $(\xi_{p_{(i)}})_x,\,i\in I_{d(n-1)}$ are a basis of $T_x(M_c(n))$ (see
(3.20)). }
\vs {\bf Proof.} The tangent vectors in the statement of the theorem are obviously a basis of $T_x(A\cdot x)$ since
$A\cdot x$ is a maximal orbit. But then the proposition follows from the equality $A\cdot x = M_c(n)$. QED\vs
The subalgebra $J(n)\s P(n)$ (see (2.30)) is Poisson commutative by (2.31). Since $P(n)$ is a Poisson algebra, 
the set of all
$f\in P(n)$ which Poisson commutes with all polynomials in $J(n)$ is an algebra containing $J(n)$. We now assert that
that this set is equal to $J(n)$.\vs {\bf Theorem 3.25.} {\it Let $J(n)\s P(n)$ be the subalgebra (and
polynomial ring in $d(n)$ generators) defined by  (2.30). Then $J(n)$ is a maximal Poisson commutative subalgeba of
$P(n)$.}\vs {\bf Proof.} Since $e +\b$ is a translate of linear subspace of $M(n)$ it follows from Theorem 2.3 that the
affine ring of $e +\b$ is the polynomial ring in the restrictions $p_i|(e +\b),\,i\in I_{d(n)}$. Thus for any $f\in
P(n)$ there exists $p\in J(n)$ such that $f = p$ on $e + \b$. But now if $f$ Poisson commutes with all elements in
$J(n)$ then $f|M_c(n)$ is a constant by Proposition 3.24 and the connectivity of $M_c(n)$, implied by Theorem 3.23.
Thus
$f = p$ on $M_{\Omega}(n)$. But $M_{\Omega}(n)$ is Zariski dense in $M(n)$. Hence $f = p$. QED\vs  We recall (see
(3.51)) $\u$ is the Lie algebra of all lower triangular matrices. For any $x\in M(n)$ let $x_{\u}$ be the ``component"
of $x$ in $\u$. That is $x_{\u}\in \u$ is such that $x-x_{\u} \in \b^T$. Noting Proposition 3.18 one has $$x_{\u} =
\sum_{m\in I_{n-1}} x_{\{m\}}\eqno (3.56)$$ Clearly (3.56) is just the decomposition of $x_{\u}$ defined by the
direct sum (3.51). Now let $c\in \Bbb C^{d(n)}$ and let $$\beta_c:M_c(n)\to \u\eqno (3.57)$$ be the regular morphism  
defined by putting $\beta_c(x) = x_{\u}$. Of course $\u$ is a nonsingular variety of dimension $d(n-1)$. If $c\in
\Omega_n$ then $M_c(n)$ is also a nonsingular variety of dimension $d(n-1)$ by Theorem 3.23. For such $c$ we can
establish the following description of $M_c(n)$. 
\vs {\bf Theorem 3.26.} {\it Let $c\in \Omega_n$ (see (2.53)). Let $\u_c\s \u$ be the image of $\beta_c$. Then $\u_c$ is
a nonempty Zariski open subset of $\u$ and $$\beta_c:M_c(n)\to \u_c\eqno (3.58)$$ is an algebraic isomorphism. }\vs
{\bf Proof. } We first prove that $\beta_c$ is injective. Let $x,y\in M_c(n)$ and assume $x_{\u} = y_{\u}$. But then
$x_{\{m\}} = y_{\{m\}}$ for all $m\in I_n$. To prove that $x=y$ we will inductively (upward) prove that $x_m = y_m$ for
all $m \in I_n$. Since $x_m$ and $y_m$ have the same spectrum one has $tr\,x_m = tr\,y_m$ for all $m\in I_n$ and hence
$\alpha_{j\,j}(x) = \alpha_{j\,j}(y)$ for all $j\in I_n$. In particular this holds for $j=1$. Thus $x_1 = y_1$. Assume
$m\in I_{n-1}$ and $x_m = y_m$. Let $g\in GL(m)$ be such that $g\,x_m g^{-1}$ is diagonal in $M(m)$. We use the
notation and arguments in the proof of Theorem 2.17. The characteristic polynomial of $x_{m+1}$ and also of $y_{m+1}$
is given by (2.59). But, from (2.59), the residue of $p(\lambda)/\prod_{i=1}^m (\lambda-\mu_{i\,m})$ at $\lambda =
\mu_{i\,m}$, where $i\in I_m$, is $-a_{i\,m+1}b_{m+1\,i}$. This is nonzero by (2.60) and is unchanged if $y$ replaces
$x$. But
$x_{\{m\}} = y_{\{m\}}$ implies that $\alpha_{i\,m+1}(g\,x_{m+1} g^{-1}) = \alpha_{i\,m+1}(g\,y_{m+1} g^{-1})$
(recalling that
$Y_m$ is stable under $Ad\,Gl(m)$). Thus $\alpha_{m+1\,i}(g\,x_{m+1} g^{-1}) = \alpha_{m+1\,i}(g\,y_{m+1} g^{-1})$.
But then $g\,x_{m+1} g^{-1} = g\,y_{m+1} g^{-1}$. Hence $x_{m+1} = y_{m+1}$. This proves the injectivity of (3.57).

We next wish to prove that the differential of (3.57) is an isomorphism at all points of $M_c(n)$. Assume not, so there
exists $x\in M_c(n),\,0\neq \xi\in \a$ such that $(\beta_c)_*(\xi_x) = 0$. Recalling (3.19) let $\xi_j$ be the
component of $\xi$ in $\a(j)$ (so that $j\in I_{n-1}$). Let $m$ be minimal such that $\xi_m\neq 0$. But now
$\xi(\alpha_{k\,m+1})(x) = 0$ for $k\in I_m$. On the other hand $\xi_j(\alpha_{k\,m+1}) = 0$ for $j\geq m+1$
since
$f_{(i,j)}\in P(j)^{Gl(j)}$ (see (2.38)) Poisson commutes with $\alpha_{k\,m+1}$ for all $i\in I_j$. Thus $\xi_m
(\alpha_{k\,m+1})(x) = 0$. But then, by (2.43) there exists $0\neq z\in Z(x,m)$ (see (2.46)) such that 
$\alpha_{k\,m+1}([z,x]) = 0$ for all $k\in I_m$. Let $W(m)$ be the $B$-orthocomplement of $Y(m)^T$ in $M(n)$ so that
one has $$M(n) = Y(m)\oplus W(m)\eqno (3.59)$$ (see (3.40)). It is clear that both summands in (3.59) are stable 
under $ad\,M(m)$. But the component of $x$ in $Y(m)$ relative to (3.59) is $x_{\{m\}}$. But $\alpha_{k\,m+1}$
clearly vanishes on $W(m)$. Thus $\alpha_{k\,m+1}([z,x_{\{m\}}]) = 0$ for all $k\in I_m$. But this implies
$[z,x_{\{m\}}]= 0$. However $x_{\{m\}}\in Y^x(m)$ by Proposition 3.18. But then $[z,x_{\{m\}}]=0$ clearly contradicts
(e) of Proposition 3.17. This proves that the differential of (3.57) is everywhere a linear isomorphism. The fact that
$\u_c$ is Zariski open in $\u$ follows by combining III, Proposition 10.4, p. 270 in [H] with III, Exercise 9.1, p.266
in [H]. (I thank P. Etingof for this reference). But then (3.58) is an isomorphism by Theorem 5.2.8, p. 85 in [S]. 
QED\vs  3.6. Let $\a^{diag}$ be the $n-1$ dimensional subalgebra spanned by $\xi_{(m,m)},\, m\in I_{n-1}$ and let
$A^{diag}\s A$ be the corresponding $n-1$ dimensional subgroup. In Remark 2.4 it was noted that the $p_j,\,j\in I_n$,
are invariant under
$Ad\,Diag(n)$. Consequently the same statement is true for any $p\in J(n)$ and hence $Ad\,Diag(n)$ necessarily
commutes with the action of $A$. We now oserve that, in fact, that $A^{diag}(n)$ operates on $M(n)$ as $Ad\,Diag(n)$ so
that $A$ can be thought of as an $d(n-1)$-dimensional extension of $Ad\,Diag(n)$. Let $\d(n)$ be the Lie algebra of
all diagonal matices in $M(n)$ so that $\d(n) = Lie \, Diag(n)$.  Actually we have already, analogously defined
$\d(m)\s M(m)$ for
$m\in I_{n-1}$. See (2.61). Clearly $I_m,\,m\in I_n,$ is a basis of $\d(n)$. But $ad\, Id_n$ operates trivially
on $M(n)$ so that $Ad\,Diag(n)$ is an $n-1$-dimensional group. Let $$\rho^{diag}: A^{diag}\to Diag(n) \eqno (3.60)$$ be
the homomorphism whose differential, $(\rho^{diag})_*$, is such that $$ (\rho^{diag})_*(\xi_{(m,m)}) = - I_m\eqno
(3.61)$$ for $m\in I_{n-1}$. One notes then that $A^{diag}\to Ad(Diag(n)$ where $a \mapsto Ad\,\, \rho^{diag}(a)$ is an
epimorphism.\vs {\bf Theorem 3.27. } {\it Let $a\in A^{diag}$ and $x\in M(n)$. Then one has $$a\cdot x =
Ad\,\, \rho^{diag}(a)(x)\eqno (3.62)$$ so that the action of $Ad\,Diag(n)$ on $M(n)$ is realized by the subgroup
$A^{diag}$ of $A$ on $M(n)$. } \vs {\bf Proof.} Write $a = a(1)\cdots a(n-1)$ where $a(m)\in A(m)$ (see (3.21)). Then
there exists $b(m)\in \Bbb C$ such that $a(m) = exp\,\,b(m)\,\xi_{(m,m)}$. Let $g(m) = \rho_{x,m}(a(m))$ (see
(3.7)). Then, by Theorem 3.6, $a\cdot x= Ad\,h(x)$ where  $h = g(1)\cdots g(m-1)$. Of course $a(m)\in A^{diag}$ and
hence it suffices to show that $\rho_{x,m}(a(m)) = \rho^{diag}(a(m))$ (which implies $g(m)$, in this case, is
independent of $x$). But to prove this it suffices to show that $$(\rho_{x,m})_*(\xi_{(m,m)}) =
(\rho^{diag})_*(\xi_{(m,m)})\eqno (3.63)$$ But indeed both sides of (3.63) are equal to $-Id_m$ by (3.8) and (3.61).
 QED\vs Let $c\in \Omega(n)$ (see (2.53)). Let $$D_c = \{a\in A\mid\,\hbox{$a$ operates as the identity map
on
$M_c(n)$}\}$$ Since $A$ is abelian and operates transitively on $M_c(n)$ (see Theorem 3.23) it follows that $$D_c=
\{a\in A\mid \hbox{ there exists $x\in M_c(n)$ such that $a\cdot x = x$}\}\eqno (3.64)$$ For $m\in I_{n-1}$ let
$D_c(m) = D_c\cap A(m)$. Let $A_c = A/D_c$. \vs {\bf Theorem 3.28.} {\it Let $c\in \Omega(n)$ (see (2.53)) and let $x\in
M_c(n)$. Then (see (3.7))
$D_c(m) = Ker\,\rho_{x,m}$. Moreover $D_c(m)$ is a closed discrete subgroup of $A(m)$ and $$\eqalign{A(m)/D_c(m)&\cong
G_{x,m}\cr &\cong
\Bbb (C^{\times})^m\cr}\eqno (3.65)$$ giving $A(m)/D_c(m)$ the structure of an abelian reductive algebraic group (i.e.
a complex torus). In addition $D_c = D_c(1)\times\cdots \times D_c(n-1)$ so that $$A_c\cong (A(1)/D_c(1)) \times 
\cdots\times
(A(n-1)/D_c(n-1))\eqno (3.66)$$ so that $A_c$ has the structure of a complex torus of dimension $d(n-1)$ which
operates simply and transitively on $M_c(n)$.} \vs {\bf Proof.} One has $Ker\,\rho_{x,m}\s D_c$ by (3.11) where we put
$y= x$. On the other hand let $a\in A(m)$ and let $g = \rho_{x,m}(a)$ so that $a\cdot x= Ad\,g(x)$ by (3.11). But if
$Ad\,g(x) = x$ then $g = 1$ by Theorem 3.14 where we put $g=g(m)$ and $g(k)=1$ for $k\neq m$ (see (3.35)). This proves
$D_c(m) = Ker\,\rho_{x,m}$, using (3.64). Of course $D_c(m)$ is closed and discrete in $A(m)$ since $(\rho_{x,m})_*$ is
an isomorphism (see (3.8)). The equalities (3.65) then follow from (3.38). Now let $a\in D_c$ so that $a\cdot x =1$.
Write $a= a(1)\cdots a(n-1)$ where $a(m)\in A(m)$. Let $g(m) = \rho_{x,m}(a(m))$. Then $g(m)=1$ for all $m\in I_{n-1}$
by Theorem 3.14 (see (3.35)). Hence $a(m)\in D_c(m)$ for all $m\in I_{n-1}$. But this readily implies (3.66). The
last statement then follows from Theorem 3.23. QED\vs Let
$c\in
\Omega(n)$ and let
$D_c^{diag} = D_c\cap A^{diag}$.\vs {\bf Theorem 3.29.} {\it Let $c\in \Omega(n)$. Then independent of $c$ 
one has (see (3.60)) $$D_c^{diag} = Ker\,\rho^{diag}\eqno (3.67)$$ and $\rho^{diag}$ induces an isomorphism
$$\eqalign{A^{diag}/ D_c^{diag}& \cong Ad\,Diag(n)\cr &\cong (\Bbb C^{\times})^{n-1}\cr}\eqno (3.68)$$ In particular
$Ad\,Diag(n)$ operates faithfully and without fixed point on $M_c(n)$. }\vs {\bf Proof.} Let $x\in M_c(n)$. Obviously
$Ker\,\rho ^{diag}\s D_c^{diag}$ by (3.62). Conversely let
$a\in D_c^{diag}$. We will use the notation and arguments in the proof of Theorem 3.27. One has $a\cdot x = x$. But then
$g(m) =1$ for all
$m\in I_{n-1}$ by Theorem 3.14. But, as established in the proof of Theorem 3.27, $g(m) = \rho^{diag}(a(m))$ Thus
$a(m)\in
 Ker\,\rho^{diag}$ for all $m\in I_{n-1}$. Hence $a\in Ker\,\rho^{diag}$. But it is immediate from the definition of
$\rho^{diag}$ that $$Diag(n) = Cent\,Gl(n)\times Im\,\rho^{diag}\eqno (3.69)$$ But this readily implies (3.68). QED\vs
3.7. We recall that the operation of transpose ($x\mapsto x^T$) in $M(n)$ stabilizes $M_c(n)$ for any $c\in 
\Bbb C^{d(n)}$. The relation of this operation to the action of $A$ is given in \vs {\bf Proposition 3.30.} {\it Let
$a\in A$ and let $x\in M(n)$. Then $$(a\cdot x)^T = a^{-1}\cdot x^T\eqno (3.70)$$}\vs {\bf Proof.} Write $a =
a(1)\cdots a(n-1)$ where $a(m)\in A(m)$ for $m\in I_{n-1}$. Let $g(m) = \rho_{x,m}(a(m))$. Then for $y\in x_m +
M(m)^{\perp}$ one has $a(m)\cdot y= g(m)\,y \,g(m)^{-1}$, by (3.11), so that $$(a(m)\cdot y)^T = (g(m)^T)^{-1}\,
y^T\,g(m)^T\eqno (3.71)$$ since $(g(m)^T)^{-1} = (g(m)^{-1})^T$. Let $g_T(m) = \rho_{x^T,m}(a(m))$. On the other hand
$(x_m)^T = (x^T)_m$ and hence $y^T \in (x^T)_m + M(m)^{\perp}$. Thus $$(a(m))^{-1}\cdot y^T =
(g_T(m))^{-1}\,y^T\,g_T(m)\eqno (3.72)$$ But we now assert that $$g(m)^T = g_T(m)\eqno (3.73)$$ Indeed, by (3.8),
$(\rho_{x^T,m})_*(\xi_{(k,m)}) = -((x^T)_m)^{m-k}$. But clearly $((x^T)_m)^{m-k}= ((x_m)^{m-k})^T$. Thus 
$$(\rho_{x^T,m})_*(\xi)= (\rho_{x,m})_*)(\xi)^T\eqno (3.74)$$ for any $\xi\in \a(m)$. But this clearly implies (3.73).
But then (3.71) and (3.72) yield
$$(a(m)\cdot y)^T = a(m)^{-1}\cdot y^T\eqno (3.75)$$ for any $y\in x_m + M(m)^{\perp}$. But since $x$ is arbitrary one
has (3.75) for all $y\in M(n)$ and all $m\in I_{n-1}$. However since $A$ is commutative this immediately implies
(3.70). QED \vs {\bf Remark 3.31.} If $c\in\Omega(n)$ the action of $A$ on $M_c(n)$ descends to an action of the
$d(n-1)$-dimensional torus $A_c$ on $M_c(n)$, by Theorem 3.28. If, in Proposition 3.30, one has $x\in M_c(n)$, where 
$c\in\Omega(n)$, one clearly has (3.70) where we can regard $a\in A_c$. \vs It is suggestive from Theorems 2.19 and
2.21 that if $x\in M(n)$ satisfies the eigenvalue disjointness condition (i.e. $x\in M_{\Omega}(n)$) and $x$ is
symmetric (i.e. $x = x^T$) then $x$ is ``essentially" uniquely determined by the set of eigenvalues $E_x(m)$ of $x_m$
for all $m\in I_n$. The precise statement is given in Theorem 3.32 below. If $c\in \Omega(n)$ let $F_c$ be
the group of elements of order $\leq 2$ in $A_c$ so that $F_c$ is a group of order $2^{d(n-1)}$ and let
$M_c^{(sym)}(n)$ be the set of symmetric matrices in $M_c(n)$. \vs {\bf Theorem 3.32.} {\it Let $c\in \Omega(n)$ (see
(2.53)). Then $M_c^{(sym)}(n)$ is a finite set of cardinality $2^{d(n-1)}$. In fact $M_c^{(sym)}(n)$ is an orbit of
$F_c$ (see Remark 3.31). } \vs {\bf Proof.} Let
$x\in M_c(n)$ so that the most general element $y\in M_c(n)$ is uniquely of the form $y = a\cdot x$ for $a\in A_c$ by
Theorem 3.28. But then the condition that $y$ be symmetric is that $a\cdot x = a^{-1}\cdot x^T$ by (3.70) and Remark
3.31. But this is just the condition that $$a^2\cdot x = x^T\eqno (3.76)$$ But now by (3.49) and Theorem 3.28 there
exists $b\in A_c$ such that $b\cdot x = x^T$. But since $A_c\cong \Bbb (C^{\times})^{d(n-1)}$ there exists $a\in A_c$
such that $b = a^2$. But then $y \in M_c^{(sym)}(n)$ where $y= a\cdot x$. Thus $M_c^{(sym)}(n)\neq \emptyset$. Now
choose $x\in M_c^{(sym)}(n)$. But then, by (3.76), a necessary and sufficient conditon that $y\in M_c^{(sym)}(n)$ is
that $a\in F_c$. QED\vs An important special case of Theorem 3.32 is when $M_c(n)$ contains a real symmetric matrix.
Let $c\in \Bbb C^{d(n)}$. The using the notation of (2.10) we will say that $c$ satisfies the eigenvalue interlacing
condition if
$\mu_{k\,m}(c)$ is real for all $m\in I_n,\,k\in I_m$ and (2.64) is satisfied
where
$\mu_{k\,m}= \mu_{k\,m}(c)$. In such a case of course $c\in \Omega(n)$. 

For any $c\in \Bbb C^{d(n)}$ and $m\in I_n$ let $p_{c,m}(\lambda)$ be the polynomial defined by putting
$$p_{c,m}(\lambda)= \prod_{k=1}^m (\lambda-\mu_{k\,m})\eqno(3.77)$$ \vskip .5pc {\bf Remark 3.33.} Assume that $c\in
\Omega(n)$ and that $\mu_{k\,m}(c)$ is real for all $m\in I_n,\,k\in I_m$. The lexicographical ordering implies 
$$\mu_{k\,m}(c)<
\mu_{k+1\,m}(c)\,\,\,\,\, \forall m\in I_n,\,\, k\in I_{m-1}\eqno (3.78)$$ Note that if $c$ satisfies the
eigenvalue interlacing condition, then for 
$m\in I_{n-1}$, $$sign(p_{c,m+1}(\mu_{i\,m}(c))) = (-1)^{m+1-i}\,\,\,\,\,\forall i\in I_{m}\eqno (3.79)$$ \vskip .5pc
{\bf Theorem 3.34.} {\it Let $c\in \Omega(n)$. Then the following conditions are equivalent
$$\eqalign{&(1)\,\,\,\hbox{$c$ satifies the eigenvalue interlacing condition}\cr &(2)\,\,\,\hbox{there exists a real
symmetric matrix in
$M_c^{(sym)}(n)$}\cr  &(3)\,\,\,\hbox{all $2^{d(n-1)}$ matrices in $M_c^{(sym)}(n)$ are real symmetric}\cr}\eqno
(3.80)$$}\vs {\bf Proof.} Obviously (3) implies (2). Assume (2) and let $x\in M_c^{(sym)}(n)$ be real symmetric. Then
$x$ is strongly regular by (2.55). Thus (2) implies (1) by Proposition 2.18. Now assume (1) and let $x\in
M_c^{(sym)}(n)$. We wish to show that
$x$ is real. But now (3.78) is satisfied. Obviously $x_1$ is real since $x_1 = \mu_{1\,1}(c)$.
Assume inductively that $m\in I_{n-1}$ and that $x_m$ is real. We use the notation and argument in the proof of Theorem
2.17. We can take $g$ to be an orthogonal matrix. Then $g\,x_{m+1} g^{-1} $ is still symmetric so that $a_{j\,m+1} = 
b_{m+1\,j}$. It suffices then to show that $a_{i\,m+1}$ is real for $i\in I_m$. (One has $\alpha_{m+1\,m+1}(x)$
is real since the reality of $E_{x}(m+1)$ implies that $tr\,x_{m+1}$ is real.) But $p(\lambda) =
p_{c,m+1}(\lambda)$. See (2.59). But then $$p_{c,m+1}(\mu_{i,m}(c)) =
-a_{i\,m+1}^2\,\prod_{j\in I_m-\{i\}}(\mu_{i\,m}-\mu_{j\,m})\eqno (3.81)$$ But clearly $$sign(-\prod_{j\in
I_m-\{i\}}(\mu_{i\,m}-\mu_{j\,m}) = (-1)^{m+1-i}$$ Thus $a_{i\,m+1}^2>0$ by (3.79). Hence $a_{i\,m+1}$ is real for all
$i\in I_m$. QED\vs Let $c\in \Omega(n)$. Recalling Theorem 3.29 put $A_c^{diag}= A^{diag}/D_c^{diag}$ so that
$A_c^{diag}$ is a $n-1$-dimensional subtorus of $A_c$, by (3.68). Now let $F_c^{diag}$ be the group of all elements
$a\in A_c^{diag}$ of order $\leq 2$ so that $F_c^{diag}$ is a subgroup of $F_c$ of order $2^{n-1}$. The determination
of $a\cdot x$ for $a\in F_c$ and $x\in M_c(n)$ seems quite nontransparent to us. However if $a\in F_c^{diag}$ then it
is easy to exhibit $a\cdot x$. Indeed it follows from Theorem 3.29 (see also (3.69)) that for each $a\in F_c^{diag}$
and $j\in I_{n-1}$ there exists $\varepsilon_j(a)\in \{1,-1\}$ such that for the faithful descent of $\rho^{diag}$
to
$F_c^{diag}$ one has $$\rho^{diag}(a) = diag(\varepsilon_1(a),\ldots,\varepsilon_{n-1}(a),1)\eqno (3.82)$$ and
$$a\cdot x = \rho^{diag}(a)\,x\, \rho^{diag}(a) \eqno (3.83)$$ \vskip .5pc {\bf Remark 3.35.} Let $x\in M(n)$ be a
real symmetric Jacobi matrix (e.g. arising say from orthogonal polynomials on $\Bbb R$- see Theorem 2.20). Let $c\in
\Bbb C^{d(n)}$ be such that $x\in M_c(n)$ so that $c$ is eigenvalue interlacing by Theorem 2.19. Then note that 
$$\{ a\cdot x\mid a\in F_c^{diag}\}\,\,\,\hbox{ is the set of all symmetric Jacobi matrices in $M_c^{(sym)}$}
\eqno (3.84)$$ by Theorem 2.19. It is interesting to note that if $a\in F_c-F_c^{diag}$ then $a\cdot x \in
M_c^{(sym)}$ is non--Jacobi.\vs {\bf Example}.  Consider the case where $n=3$ so that there are
8 symmetric matrices in $M_c(3)$ for any $c\in \Omega(n)$. Let $c$ be defined so that $E_c(1) = \{0\},\,E_c(2)
= \{1,-1\},\,E_c(3) = \{\sqrt 2,0,-\sqrt 2\}$ so that $c$ is eigenvalue interlacing.
Then 2 of the 8 real symmetric matrices in $M_c(3)$ are 
$$x=\left(\matrix{0&1&0\cr 1&0&1\cr 0&1&0\cr}\right),\qquad
y=\left(\matrix{0&1&1\cr 1&0&0\cr 1&0&0\cr}\right)$$ noting that $x$, but not $y$, is Jacobi. The remaining 6 are
obtained by sign changes in $x$ and $y$. That is, by applying the three non-trivial elements in $F_c^{diag}$ to $x$ and
$y$.\vs 3.8. In Part II we will need the following result on polarizations of regular adjoint orbits and maximal
$A$-orbits. Let $x\in M(n)$. It is clear from (1.18) that the adjoint $O_x$ is stable under the action of $A$. Let
$O_x^{sreg} = M^{sreg}\cap O_x$ so that $O_x^{sreg}$ is Zariski open subset of $O_x$.\vs {\bf Theorem 3.36.} {\it Let
$x\in M(n)$. Then $O_x^{sreg}$ is non-empty (and hence dense in $O_x$) if and only if $x$ is regular in $M(n)$. In
particular if
$x$ is regular then $O_x^{sreg}$ is a symplectic $2\,d(n-1)$-dimensional manifold (in the complex sense). Furthermore
$O_x^{sreg}$ is stable under the action of $A$ and the orbits of $A$ in $O_x^{sreg}$ (necessarily of
dimension $d(n-1)$) are the leaves of a polarization of $O_x^{sreg}$.}\vs {\bf Proof.} If $x$ is not regular then all
elements in
$O_x$ are not regular so that $O_x^{sreg}$ is empty (strongly regular implies regular -- see Theorem 2.13). Now assume
that
$x$ is regular. But now by Theorem 2.5 there exists $y\in e + \b$ such that $E_x(n) = E_y(n)$.
But $y$ is regular by (2.27). But, as one knows, any two regular matrices with the same spectrum are conjugate (they
are both conjugate to the same companion matrix). Thus $y\in O_x$. But then $y\in O_x^{sreg}$ by (2.27). Hence $
O_x^{sreg}$ is not empty. But then $O_x^{sreg}$ is a union of maximal $A$-orbits by (3.30). But each such orbit is a
Lagrangian submanifold of $O_x$ by dimension (see (3.29)) and the vanishing of (1.17), by (2.31), where $\psi =
p_{(i)}$ and $u= (\xi_{p_{(j)}})_x$ (see (1.18)) for $i,j\in I_{d(n-1)}$. QED\vs

\centerline{\bf References}\vskip 1pc
\rm
\item {[B]} A. Borel, {\it Linear Algebraic Groups}, W. A. Benjamin,
Inc, 1969
\item {[C]} C. Chevalley, {\it Fondements de la G\'eom\'etrie
Alg\'ebrique}, Facult\'e des Sciences de Paris, Math\'ematiques
approfondies, 1957/1958
\item {[CG]} N. Chris and V.Ginzburg, {\it Representation Theory
and Complex Geometry}, Birkh\"auser, 1997
\item {[GS]} V. Guillemin and S. Sternberg, On the collective
integrability according to the method of Thimm, Ergod. Th. \& Dynam.
Sys. (1983), {\bf 3}, 219--230
\item {[H]} R. Hartshorne, {\it Algebraic Geometry}, Grad. Texts in
Math., Vol.  52, Springer-Verlag, 1977
\item {[J]} D. Jackson, {\it Fourier Series and Orthogonal
Polynomials},  The Carus Mathematical Monographs, No. 6, Math.
Assoc. of America, 1948
\item {[K]} B. Kostant, Quantization and Unitary Representations,
Lecture Notes in Math, Vol.  170, Springer-Verlag, 1970
\item {[K-1]} B. Kostant, Lie group representations on
polynomial rings, AJM, Vol. 85(1963), 327-404
\item {[M]} D. Mumford, {\it The red book of varieties and
schemes}, Lecture Notes in Math., Vol. 1358, Springer, 1995
\item {[S]} T. Springer, {\it Linear Algebraic Groups}, 2nd
edition, Prog. Math., Vol.  9, Birkh\"auser Boston, 1998
\parindent=30pt
\baselineskip=14pt
\vskip 1pc
$$\halign{\hskip 8pc # \qquad\hfill & #\hfill & #\hfill & #\hfill & # \qquad \hfill\cr
Bertram Kostant & Nolan Wallach\cr
Dept. of Math. & Dept. of Math.\cr
MIT & UCSD \cr
Cambridge, MA 02139 & San Diego, CA 92093 \cr
kostant@math.mit.edu & nwallach@ucsd.edu\cr}
$$
\end

\vbox to 60pt{\hbox{Bertram Kostant}
      \hbox{Dept. of Math.}
      \hbox{MIT}
      \hbox{Cambridge, MA 02139}} \noindent E-mail
kostant@math.mit.edu\vskip 1pc
      
\vs\vbox to 60pt{\hbox{Nolan Wallach}
      \hbox{Dept. of Math.}
      \hbox{UCSD}
      \hbox{San Diego, CA 92093}}\noindent
       E-mail nwallach@ucsd.edu

\end